\documentclass[a4paper, reqno, 11pt]{amsart}
\usepackage{amscd, amssymb, amsthm, amsmath}
\usepackage{graphicx, psfrag, diagbox}
\newtheorem*{MainA}{Main Theorem A}
\newtheorem*{MainB}{Main Theorem B}

\newtheorem{lem}{Lemma}[section]
\newtheorem{thm}[lem]{Theorem}

\newtheorem{cor}[lem]{Corollary}

\newtheorem{defi}[lem]{Definition}

\newcommand{\Stat}{{\textsf{Stat}}}
\newcommand{\stat}{{\textsf{stat}}}
\newcommand{\maj}{{\textsf{maj}}}
\newcommand{\inv}{{\textsf{inv}}}
\newcommand{\Nega}{{\textsf{Neg}}}
\newcommand{\nega}{{\textsf{neg}}}

\newcommand{\fmaj}{{\textsf{fmaj}}}
\newcommand{\Fmaj}{{\textsf{Fmaj}}}
\newcommand{\nmaj}{{\textsf{nmaj}}}
\newcommand{\Dmaj}{{\textsf{Dmaj}}}
\newcommand{\dmaj}{{\textsf{dmaj}}}
\newcommand{\lmaj}{{\textsf{lmaj}}}
\newcommand{\rmaj}{{\textsf{rmaj}}}
\newcommand{\fmaf}{{\textsf{fmaf}}}
\newcommand{\Fix}{{\textsf{Fix}}}
\newcommand{\fix}{{\textsf{fix}}}
\newcommand{\rinv}{{\textsf{rinv}}}
\newcommand{\col}{{\textsf{col}}}
\newcommand{\sor}{{\textsf{sor}}}

\numberwithin{equation}{section}

\renewcommand{\d}{\displaystyle}
\def\proof{\noindent{\bf Proof.\quad}}
\def\qed{\hspace*{\fill} $\Box$ \medskip}

\begin{document}
%%%%%%%%%%%%%%%%%%%%%%%%%%%%%%%%%%%%%%%%%%%%%%%%%%%%%%%%%%%
\title[Signed Mahonian on $G(r,1,n)$]{Signed Mahonian polynomials for major and sorting indices}

\author{Huilan Chang}
\address{Department of Applied Mathematics \\
National University of Kaohsiung \\
Taiwan, ROC} 
\email[Huilan Chang]{huilan0102@gmail.com}

\author{Sen-Peng Eu}
\address{Department of Mathematics \\
National Taiwan Normal University \\
Taiwan, ROC} 
\email[Sen-Peng Eu]{speu@ntnu.edu.tw}

\author{Shishuo Fu}
\address{College of Mathematics and Statistics \\
Chongqing University \\
Chongqing, P.~R.~China}
\email[Shishuo Fu]{fsshuo@cqu.edu.cn}

\author{Zhicong Lin}
\address{School of Science \\
Jimei University \\
Xiamen, P.~R.~China}
\email[Zhicong Lin]{zhicong.lin@univie.ac.at}

\author{Yuan-Hsun Lo}
\address{School of Mathematical Sciences \\
Xiamen University \\
Xiamen, P.~R.~China}
\email[Yuan-Hsun Lo]{yhlo0830@gmail.com}



\thanks{Partially supported by Ministry of Science and Technology, Taiwan under grants 101-2115-M-390-004-MY3 (S.-P.~Eu) and 106-2115-M-390-004-MY2 (H.~Chang), the National Natural Science Foundation of China under grant numbers 11501244 (Z.~Lin), 11501061 (S.~Fu), 11601454 (Y.-H.~Lo) and 11871247 (Z.~Lin), the Natural Science Foundation of Fujian Province, China under grant number 2016J05021 (Y.-H.~Lo), and the Fundamental Research Funds for the Central Universities in China under grant number 20720150210 (Y.-H.~Lo)}
%\date{\today}
\maketitle
%%%%%%%%%%%%%%%%%%%%%%%%%%%%%%%%%%%%%%%%%%%%%%%%%%%%%%%%%%%

\begin{abstract}
%We characterize the 1-dim characters of the complex reflection group $G(r,1,n)$ in terms of the length function.
We derive some new signed Mahonian polynomials over the complex reflection group $G(r,1,n)=C_r\wr\mathfrak{S}_n$, where the ``sign'' is taken to be any of the $2r$ $1$-dim characters and the ``Mahonian'' statistics are the $\lmaj$ defined by Bagno and the $\sor$ defined by Eu \emph{et al}. 
Various new signed Mahonian polynomials over Coxeter groups of types $B_n$ and $D_n$ are obtained as well.
We also investigate the signed counting polynomials on $G(r,1,n)$ for those statistics with the distribution $[r]_q[2r]_q\cdots [nr]_q$.
\end{abstract}

%%%%%%%%%%%%%%%%%%%%%%%%%%%%%%%%%%%%%%%%%%%%%%%%%%%%%%%%%%%%%%%%%%%%%%%%%%

\section{Introduction}\label{sec:intro}

%%%%%%%%%%%%
\subsection{Signed Mahonian}

Let $\mathfrak{S}_n$ be the symmetric group of $\{1,2,\dots , n\}$. The
\emph{inversion} and $\emph{major}$ statistics of a permutation $\pi
= \pi_1\pi_2\dots \pi_n \in \mathfrak{S}_n$ are defined respectively by
\begin{align*}
\inv(\pi):&=| \{(i,j): i<j \mbox{ and } \pi_i>\pi_j \}|,\\
\maj(\pi):&= \sum_{\pi_i>\pi_{i+1}} i .
\end{align*}
A fundamental result of MacMahon~\cite{MacMahon_13} states that
$\inv$ and $\maj$ have the same distribution over $\mathfrak{S}_n$ together
with the generating function
\begin{equation}\label{eq:A_Poincare}
\sum_{\pi\in \mathfrak{S}_n}q^{\inv(\pi)}= \sum_{\pi\in \mathfrak{S}_n}q^{\maj(\pi)}=[1]_q[2]_q\dots [n]_q,
\end{equation}
where $[k]_q:=1+q+q^2+\dots +q^{k-1}$ is the $q$-analogue of $k$. To
celebrate this result, any statistic equidistributed with $\inv$ on
$\mathfrak{S}_n$ is called \emph{Mahonian}.

The motivation of this paper is the following ``signed Mahonian" identity \cite{Wachs_92}, obtained by Gessel and Simion.
\begin{equation} \label{eq:A_signMaho}
\sum_{\pi\in \mathfrak{S}_n}(-1)^{\inv(\pi)}q^{\maj(\pi)} = [1]_q[2]_{-q}[3]_q[4]_{-q}\cdots [n]_{(-1)^{n-1} q}.
\end{equation}
%This elegant result was first proved by Gessel and Simion~\cite{Wachs_92}.

%%%%%%%%%%%%%%%%%%%
\subsection{Extensions to Coxeter groups}
One can regard $\mathfrak{S}_n$ as the Coxeter group of type $A_{n-1}$. 
For a Coxeter group, the length function $\ell(\pi)$ is the minimum number of generators required to express a group element $\pi$. 
Note that in the type $A_{n-1}$ case the length function $\ell_A(\pi)$ is exactly $\inv(\pi)$, where the subscript $A$ is to emphasize the type.
Therefore in a Coxeter group or a complex reflection group $G$ we call a statistic \emph{Mahonian} if it is equidistributed with the length function.

On the other hand, (\ref{eq:A_Poincare}) and (\ref{eq:A_signMaho})
can be seen as
$$\sum_{\pi\in \mathfrak{S}_n} \chi(\pi) q^{\maj(\pi)},$$
where $\chi$ is one of the two 1-dim characters of $\mathfrak{S}_n$, namely $\chi=1$ or $\chi(\pi)=(-1)^{\ell_A(\pi)}$. 
Therefore from this point of view one may consider the polynomial
\begin{equation}
\sum_{\pi\in G} \chi({\pi}) q^{\stat({\pi})},
\end{equation}
where $\chi$ is any 1-dim character of $G$ and $\stat$ is some Mahonian statistic.

There have been many works extending (\ref{eq:A_signMaho}) to other Coxeter groups or reflection groups. 
In the case of type $B_n$, Adin and Roichman~\cite{Adin_Roichman_01} defined the Mahonian index so called flag major, $\fmaj$, and Adin, Gessel and Roichman~\cite{Adin_Gessel_Roichman_05} obtained the signed Mahonian polynomials with respect to $\fmaj$ together with any of the four 1-dim characters ($1, (-1)^{\ell_B(\pi)}, (-1)^{\nega(\pi)}, (-1)^{\inv(|\pi|)}$) of $B_n$.
Here, $\ell_B(\pi)$ refers to the length function of $\pi=\pi_1\pi_2\cdots\pi_n\in B_n$ which has the following combinatorial interpretation
\begin{equation}\label{eq:length_Bn}
\ell_B(\pi)=\inv(\pi)-\sum_{i\in\Nega(\pi)}\pi_i,
\end{equation}
where $\Nega(\pi):=\{i: \pi_i<0\}$, $\nega(\pi)$ is the cardinality of $\Nega(\pi)$, and $|\pi|:=|\pi_1||\pi_2|\cdots|\pi_n|\in \mathfrak{S}_n$.
Meanwhile, Fire~\cite{Fire_08} obtained another set of identities by considering the Mahonian statistics $\Fmaj$ and $\nmaj$ together with $\chi=(-1)^{\ell_B(\pi)}$.

For the case of type $D_n$, there are several Mahonian statistics in literature. 
For example, $\dmaj, \Dmaj$ and flag major index $\fmaj$. 
There are only two 1-dim characters ($1$ and $(-1)^{\ell_D(\pi)}$) on $D_n$ and the signed Mahonian problems have been studied by Biagioli~\cite{Biagioli_06} in the case of $\Dmaj$ and by Biagioli
and Caselli~\cite{Biagioli_Caselli_12} in the case of $\fmaj$.

The definitions of relative statistics for types $B_n$ and $D_n$ are given in Section~\ref{sec:Bn} and Section~\ref{sec:Dn}, respectively.

%%%%%%%%%%%%%%%%%%%%%
\subsection{Extensions to complex reflection group}

One highlight of this paper is to derive new signed Mahonian polynomials on the complex reflection group $G(r,1,n)$ (denoted by $G(r,n)$ for short throughout the paper), which is isomorphic to the wreath product $C_r\wr\mathfrak{S}_n$, where $C_r$ is the cyclic group of order $r$. 
$G(r,n)$ is also called the group of $r$-colored permutations.
Note that $G(1,n)=\mathfrak{S}_n$ and $G(2,n)=B_n$.
We will derive the signed Mahonian polynomials by taking $\chi=\chi_{a,b}$ to be any of the $2r$ 1-dim characters and $\stat$ the length function $\ell$ or the Mahonian statistic $\lmaj$ defined by Bagno~\cite{Bagno_04}.

%\begin{MainA}[Theorem~\ref{thm:1-dim}]
%$G(r,n)$ has $2r$ $1$-dim characters
%$$\chi_{a,b}(\pi)=(-1)^{a(\ell(\pi)-\col(\pi))}\omega^{b\,\col(\pi)},$$
%where $\omega$ is a primitive $r$th root of $1$, $\col(\pi)$ is the sum of colors of $\pi$, $a=0,1$ and $b=0,1,\cdots, r-1$.
%\end{MainA}

%Take $r=2$ for an example.
%In this case $\omega=-1$, $\col(\pi)=\nega(\pi)$ (the number of negatives of $\pi$), and $\ell(\pi)=\ell_B(\pi)$ by \eqref{eq:length_Bn} and \eqref{eq:length_Grn}.
%Then 
%\begin{align*}
%\chi_{0,0}(\pi) &= 1, \\
%\chi_{0,1}(\pi) &= (-1)^{\nega(\pi)}, \\
%\chi_{1,1}(\pi) &= (-1)^{\ell_B(\pi)}, \\
%\chi_{1,0}(\pi) &= (-1)^{\ell_B(\pi)-\nega(\pi)} = (-1)^{\inv(\pi)-\big(\sum_{i\in\Nega(\pi)}\pi_i\big) - \nega(\pi)} \\
%&= (-1)^{\inv(\pi)}(-1)^{\sum_{i\in\Nega(\pi)}(|\pi_i|-1)} \stackrel{(*)}{=} (-1)^{\inv(\pi)}(-1)^{\inv(|\pi|)+\inv(\pi)} \\
%&= (-1)^{\inv(|\pi|)},
%\end{align*}
%where $(*)$ is due to the identity \eqref{eq:proof_Bn_1} in the proof of Lemma~\ref{lem:inv}.
%We can see these four characters meets that in \cite{Adin_Gessel_Roichman_05}.

%By means of the description of 1-dim characters as stated in Equation~\eqref{eq:1-dim}, we will derive the signed Mahonian polynomials by taking $\chi=\chi_{a,b}$ to be any of the $2r$ $1$-dim characters and $\stat$ the length function $\ell$ or the Mahonian statistic $\lmaj$ defined by Bagno~\cite{Bagno_04}.

\begin{MainA}[Theorem~\ref{thm:G_length}]
For $a=0,1$ and $b=0,1,\cdots, r-1$, we have the following signed Mahonian polynomials:
$$\sum_{\pi\in G(r,n)} \chi_{a,b}(\pi)q^{\ell(\pi)}=\prod_{k=1}^n[k]_{(-1)^aq}\left(1+(-1)^{a(k-1)}\omega^bq^k[r-1]_{\omega^bq}\right),$$
and
$$\sum_{\pi\in G(r,n)} \chi_{a,b}(\pi)q^{\lmaj(\pi)}=\prod_{k=1}^n[k]_{(-1)^{a(k-1)}q}\left(1+(-1)^{a(k-1)}\omega^bq^k[r-1]_{\omega^bq}\right),$$
where $\omega$ is a primitive $r$th root of unity.
\end{MainA}

\begin{figure}[h]
\centering
\includegraphics[width=3.3in]{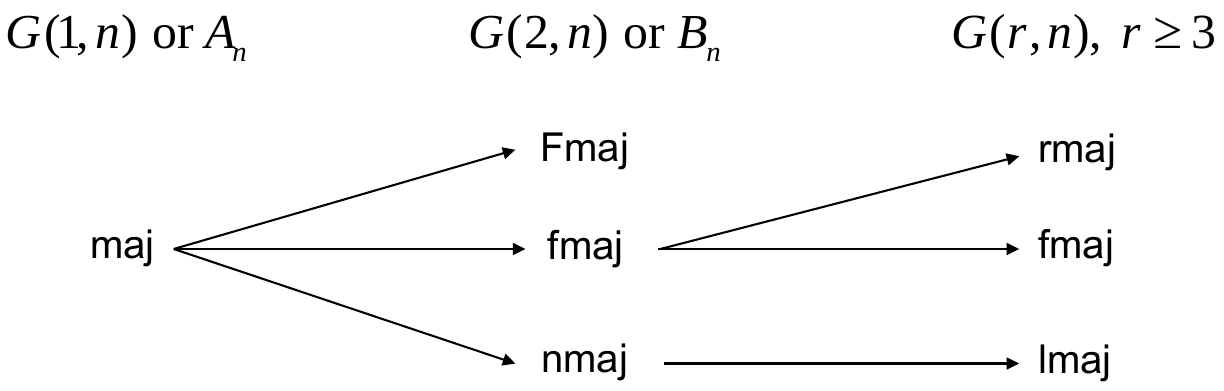}
\caption{The relationship between various major indices.} \label{fig:major}
\end{figure}

It is worth mentioning that Biagioli and Caselli~\cite{Biagioli_Caselli_12} obtained a closed form of the polynomial $\sum_{\pi\in G(r,n)}\chi(\pi)q^{\fmaj(\pi)}$ for any $1$-dim character $\chi$ of $G(r,n)$ and the flag major index $\fmaj$ defined by Adin and Roichman~\cite{Adin_Roichman_01}.
However this $\fmaj$ is not Mahonian if $r\ge 3$ (i.e., is not equidistributed with $\ell$), hence it does not equal $\lmaj$ and therefore the result does not overlap ours.
Another interesting statistic $\rmaj$, called \emph{root major index}, on $G(r,n)$ defined by Haglund, Loehr and Remmel \cite{Haglund_Loehr_Remmel_05} is equidistributed with $\fmaj$, and then is not Mahonian.
Hence $\rmaj$ is not taken into consideration in this paper.
See Figure~\ref{fig:major} for the relationship between these ``major'' statistics on $G(r,n)$.
$\stat_1\longrightarrow\stat_2$ means $\stat_1$ is extended to $\stat_2$.
For example, $\nmaj$ is exactly $\lmaj$ when $r=2$.

%%%%%%%%%%%%%%%%%%
%\subsection{Statistics with distribution $[r]_q[2r]_q\cdots [nr]_q$}
%In addition to Mahonian statistics, there are several known statistics on $G(r,n)$ having the distribution 
%$$[r]_q[2r]_q[3r]_q\dots [nr]_q.$$
%For example, $\fmaj, \fmaf, \rmaj$ and $\rinv$ are some of them (definitions will be given later). 
%In the same spirit one can formally consider the signed counting polynomial
%$$\sum_{\pi\in G(r,n)} (-1)^{\stat_1(\pi)} q^{\stat_2(\pi)}$$
%by taking any two of these four statistics. 
%It turns out that except for two cases, which seem to have no factorial-type product formulas, we do have nice closed forms.
%
%As for the proofs, the technic used is rather simple. 
%Roughly speaking, most of the proofs are done by separating `colors' and `numbers' and applying the results in type $A_{n-1}$ or $B_n$.
%

Another interesting Mohonian statistic is the \emph{sorting index} $\sor$, which is implicitly given for type $A_{n-1}$ in a bijective proof on $\mathfrak{S}_n$ by Foata and Han~\cite{Foata_Han_09}, defined for Coxeter groups by Petersen~\cite{Petersen_11}, and then generalized by Eu \emph{et al.}~\cite{ELW_15} for complex reflection groups.
Our second class of signed Mahonian polynomials takes $\sor$ into consideration.
The definition of $\sor$ will be given later in Section~\ref{sec:Grn_sor}.

\begin{MainB}[Theorem~\ref{thm:sorting}]
For $a=0,1$ and $b=0,1,\cdots, r-1$, we have the following signed Mahonian polynomials:
\begin{align*}
\sum_{\pi\in G(r,n)} & \chi_{a,b}(\pi)q^{\sor(\pi)} \\ 
& = \prod_{k=1}^n \left( 1 + (-1)^a[k-1]_q\left( q + \omega^{br}q^k[r-1]_q \right) + \omega^b q^{2k-1}[r-1]_{\omega^b q} \right),
\end{align*}
where $\omega$ is a primitive $r$th root of unity.
\end{MainB}

%%%%%%%%%%%
\subsection{Statistics with distribution $[r]_q[2r]_q\cdots [nr]_q$}
In addition to Mahonian statistics, there are several known
statistics on $G(r,n)$ having the distribution
$$[r]_q[2r]_q[3r]_q\dots [nr]_q.$$
For example, $\fmaj, \fmaf, \rmaj$ and $\rinv$ are some of them
(definitions will be given later). In the same spirit one can
formally consider the signed counting polynomial
$$\sum_{\pi\in G(r,n)} (-1)^{\stat_1(\pi)} q^{\stat_2(\pi)}$$
by taking any two of these four statistics. It turns out that except
for two cases we do have nice closed forms. 
See Table~\ref{tab:sign_counting} for a brief summary.

\begin{table}[h]
\begin{tabular}{|c||c|c|c|c|}\hline
%\backslashbox{$\stat_1$}{$\stat_2$}& \makebox[3em]{$\fmaj$}& \makebox[3em]{$\fmaf$}& \makebox[3em]{$\rmaj$}& \makebox[3em]{$\rinv$}\\
\diagbox{$\stat_1$}{$\stat_2$}& {$\fmaj$}& {$\fmaf$}& {$\rmaj$}& {$\rinv$}\\
\hline
$\fmaj$ & -- & $\ast$Theorem~\ref{thm:fmaffmaj} & Theorem~\ref{thm:fmajrmaj} & Theorem~\ref{thm:fmajrinv} \\ \hline
$\fmaf$ & $\ast$Theorem~\ref{thm:fmaffmaj} & -- & $\ast$Theorem~\ref{thm:fmafrmaj} & $\ast$Theorem~\ref{thm:fmafrmaj} \\ \hline
$\rmaj$ & Theorem~\ref{thm:fmajrmaj} & not known & -- & Theorem~\ref{thm:rmajrinv} \\ \hline
$\rinv$ & Theorem~\ref{thm:fmajrinv} & not known & Theorem~\ref{thm:rmajrinv} & -- \\ \hline
\end{tabular}
\caption{All derived signed counting polynomials, where $\ast$ indicates that the closed forms are valid for even $r$.}
\label{tab:sign_counting}
\end{table}

The rest of the paper is organized as follows. 
The proofs of the two main results (Main Theorem A -- B) are given in Section~\ref{sec:Grn}.
Some new signed Mahonian results of type $B_n$ and $D_n$ are collected in Section~\ref{sec:Bn} and Section~\ref{sec:Dn}, respectively.
Section~\ref{sec:others} is devoted to those statistics on $G(r,n)$ with $[r]_q[2r]_q\dots [nr]_q$ as the distributions.
Finally, we propose a brief conclusion in Section~\ref{sec:summary}.

%%%%%%%%%%%%%%%%%%%%%%%%%%%%%%%%%%%%%%%%%%%%%%%%%%%%%%%%%%%%%%%%%%%%%%%%%%%%%%%
%%%%%%%%%%%%%%%%%%%%%%%%%%%%%%%%%%%%%%%%%%%%%%%%%%%%%%%%%%%%%%%%%%%%%%%%%%%%%%%
\section{Signed Mahonian on $G(r,n)$}\label{sec:Grn}

%%%%%%%%%% Adin's Wreath Product %%%%%%%%%%%%%%%%%%%%%%%%%%%%%%%%%%%%%%
\subsection{Colored permutation group}\label{sec:Grn_Grn}

Let $r,n$ be positive integers. 
The group of \emph{colored permutations} of $n$ digits with $r$ colors, denoted by $G(r,n)$, is the wreath product $C_r\wr \mathfrak{S}_n$ of the cyclic group $C_r$ with $\mathfrak{S}_n$. 
It is useful to describe the elements in two ways.
First, $G(r,n)$ is generated by the set of generators $\mathcal{S}=\{s_0,s_1,\ldots,s_{n-1}\}$, which satisfies the following relations:
\begin{align} \label{eq:wreath_product}
s_0^r &= 1, \notag \\
s_i^2 &= 1, \text{ } i=1,2,\ldots,n-1, \notag \\
(s_is_j)^2 &= 1, \text{ } |i-j|>1, \notag \\
(s_is_{i+1})^3 &= 1, \text{ } i=1,2,\ldots,n-2, \\
(s_0s_1)^{2r} &= 1. \notag
\end{align}
In fact, one can realize the generator $s_i$ ($i\neq 0$) as the transposition $(i,i+1)$ and $s_0$ as the mapping
$$ s_0(i)=\left\{ \begin{array}{ll}
\omega\cdot 1, & \text{if }i=1,\\
i, & \text{otherwise,}\\
\end{array}\right. $$
where $\omega$ is a primitive $r$th root of unity.
We can also write an element $\pi\in G(r,n)$ as an ordered pair $(z,\sigma)$ or a sequence $\sigma_1^{[z_1]}\,\sigma_2^{[z_2]}\cdots\sigma_n^{[z_n]}$ with the convention $\sigma_i^{[0]}=\sigma_i$, where $z=(z_1,z_2,\ldots,z_n)$ is an $n$-tuple of integers in $\mathbb{Z}_r\cong C_r$ and $\sigma\in \mathfrak{S}_n$, such that $\sigma_i=|\pi_i|$ and ${\pi_i}/{|\pi_i|}=\omega^{z_i}$. 
For example, $\omega 2\,\,1\,\,5\,\,\omega^3 4\,\,\omega^2 3\in G(4,5)$ can be written as $((1,0,0,3,2),21543)$ or $2^{[1]}\,1\,5\,4^{[3]}\,3^{[2]}$. 
Note that $G(1,n)=\mathfrak{S}_n$ and $G(2,n)=B_n$.
Readers are referred to \cite{Adin_Roichman_01, Bagno_04, Haglund_Loehr_Remmel_05} for more
information.

Let $\ell(\pi)$ be the length of $\pi\in G(r,n)$ with respect to $\mathcal{S}$. 
Bagno \cite{Bagno_04} gave the following combinatorial interpretation of $\ell(\pi)$:
\begin{equation}\label{eq:length_Grn}
\ell(\pi)=\inv(\pi)+\sum_{z_i>0}(|\pi_i|+z_i-1),
\end{equation}
where $\inv(\pi):=|\{(i,j):\,i<j\text{ and }\pi_i>\pi_j\}|$ with respect to the linear order
\begin{equation}\label{eq:order_wreath}
(n^{[r-1]}<\cdots<n^{[1]}) <\cdots< (1^{[r-1]}<\cdots<1^{[1]}) < (1<\cdots<n).
\end{equation}
%Note that in this section all calculations (including $\maj_A$) are with respect to the above linear order.
%The $\fmaj$ can also be defined on $G(r,n)$. 
%However it is no longer Mahonian when $r\ge 3$~\cite{Adin_Roichman_01}. 
By defining $\maj(\pi):=\sum_{\pi_i>\pi_{i+1}}i$ with respect to the above linear order, a Mahonian statistic $\lmaj$ on $G(r,n)$ was found by Bagno~\cite{Bagno_04}, defined by
\begin{equation}\label{eq:lmaj_Grn}
\lmaj(\pi):=\maj(\pi)+\sum_{z_i>0}(|\pi_i|+z_i-1).
\end{equation}
It has the generating function
\begin{equation}\label{eq:equidist_Grn}
\sum_{\pi\in G(r,n)}q^{\ell(\pi)}=\sum_{\pi\in G(r,n)}q^{\lmaj(\pi)}= \prod_{k=1}^n[k]_q(1+q^k[r-1]_q).
\end{equation}

%%%%%%%%%%%%%%%%%%%%%%%%%%
\subsection{Character and signed Mahonian of $\ell$ and $\lmaj$}\label{sec:Grn_SM}
%For $\pi=(z,\sigma)\in G(r,n)$ let $\col(\pi):=\sum_{i=1}^n z_i$.
%We also denote $|\pi|=\sigma$.
%The 1-dim characters of $G(r,n)$ have been characterized in~\cite{Biagioli_Caselli_12} as the form:
%\begin{equation}\label{eq:r_1-dim-BC}
%\chi_{a,b}(\pi)=(-1)^{a\,\inv(|\pi|)}\omega^{b\,\col(\pi)},
%\end{equation}
%where $a=0,1$ and $b=0,1,\cdots, r-1$.

First we characterize all 1-dim characters of $G(r,n)$ in terms of the length function. 
For $\pi=(z,\sigma)\in G(r,n)$ let $\col(\pi):=\sum_{i=1}^n z_i$.

%----one-dim character--------------------------------------------
\begin{lem}\label{lem:1-dim}
$G(r,n)$ has $2r$ $1$-dim characters, which can be expressed as
\begin{equation}\label{eq:r_1-dim}
\chi_{a,b}(\pi)=(-1)^{a(\ell(\pi)-\col(\pi))}\omega^{b\,\col(\pi)},
\end{equation}
for $a=0,1$ and $b=0,1,\cdots, r-1$.
\end{lem}
\proof 
Let $\chi$ be any $1$-dim character. 
We first look at the values $\chi(s_i)$, $i=1,2,\ldots,n-1$.
Since, by \eqref{eq:wreath_product}, $s_i s_{i+1} s_i = s_{i+1} s_i s_{i+1}$, $\chi(s_i)$ and $\chi(s_{i+1})$ should be the same; therefore, $\chi(s_1)=\chi(s_2)=\cdots =\chi(s_{n-1})$. 
Moreover, $s_i$ is an involution, then $\chi(s_i)=(-1)^a$ for some $a=0,1$. 
Since $s_0$ is of order $r$, $\chi(s_0)=\omega^b$ for some $b=0,1,\cdots, r-1$. 
Then there are $2r$ non-isomorphic choices for $\chi$, which are denoted by $\chi_{a,b}$ for $a=0,1$ and $b=0,1,\cdots, r-1$. 
By the definitions of $\mathcal{S}$ and $\ell$, the number of occurrences of $s_0$ needed to express $\pi$ is exactly $\sum_{i=1}^n z_i=\col(\pi)$. 
Hence the expression of $\chi_{a,b}$ in \eqref{eq:r_1-dim} fits the requirements. 
It is easy to verify that each $\chi_{a,b}$ is indeed a group homomorphism. 
This completes the proof. \qed

It should be emphasized that the characterization of the 1-dim characters of $G(r,n)$ is not original, but has been studied in~\cite{Biagioli_Caselli_12} as the form:
\begin{equation}\label{eq:r_1-dim-BC}
\chi_{a,b}(\pi)=(-1)^{a\,\inv(|\pi|)}\omega^{b\,\col(\pi)},
\end{equation}
where $|\pi|=\sigma$.
The two expressions in \eqref{eq:r_1-dim} and \eqref{eq:r_1-dim-BC} are similar, however, the former one is more helpful to derive our two main results.
In fact, by the definition of $\ell$, \eqref{eq:r_1-dim} can be obtained directly from \eqref{eq:r_1-dim-BC} by showing
\begin{equation}\label{eq:r_1-dim-connection}
(-1)^{\inv(|\pi|)} = (-1)^{\inv(\pi)+\sum_{z_i>0}(|\pi_i|-1)}.
\end{equation}
The proof of the identity above is similar to the one of Lemma~\ref{lem:inv} and is omitted here.

%---------------preparation-------------------------------
To prove our signed Mahonian results we need some more preparations.
The key idea is to uniquely decompose a colored permutation $\pi=\sigma_1^{[z_1]}\,\sigma_2^{[z_2]}\cdots\sigma_n^{[z_n]}\in G(r,n)$ as $\pi=\tau\rho$, a product of a linear ordered $\tau$ of $\sigma_1^{[z_1]},\sigma_2^{[z_2]},\ldots,\sigma_n^{[z_n]}$ and a permutation $\rho\in \mathfrak{S}_n$, where $\rho_i=j$ if and only if $\pi_i$ is the $j$th smallest element among $\{\sigma_1^{[z_1]},\sigma_2^{[z_2]},\ldots,\sigma_n^{[z_n]}\}$.
For example, $$2\,4^{[2]}\,1^{[1]}\,3^{[2]}=(4^{[2]}\,3^{[2]}\,1^{[1]}\,2)~~(4\,1\,3\,2).$$
In general, let $U_{r,n}:=\{\tau\in G(r,n):\,\tau_1<\cdots<\tau_n\}$ and then $G(r,n)=U_{r,n}\cdot \mathfrak{S}_n$.
Here we adopt the linear order \eqref{eq:order_wreath}.
It is obvious that $\col(\pi)=\col(\tau)$ if $\pi=\tau\rho\in U_{r,n}\cdot \mathfrak{S}_n$. 
Note that we also have $\inv(\pi)=\inv(\rho)$ and $\maj(\pi)=\maj(\rho)$. 
The bivariate generating function
$$F(t,q):=\sum_{\tau\in U_{r,n}} t^{\sum_{z_i>0}(|\tau_i|-1)}q^{\col(\tau)}$$
is essential and we have the following formula.

\begin{lem}\label{lem:r_bivariate} We have
$$F(t,q)=\prod_{k=1}^{n}(1+t^{k-1}q+t^{k-1}q^2+\cdots+t^{k-1}q^{r-1})=\prod_{k=1}^{n}(1+t^{k-1}q[r-1]_q).$$
\end{lem}
\proof
For $1\leq k\leq n$ we define elements $\psi_k$ of the group algebra of $G(r,n)$ by
$$\psi_k=1+s_{k-1}\cdots s_1s_0+s_{k-1}\cdots s_1(s_0)^2+\cdots+s_{k-1}\cdots s_1(s_0)^{r-1}.$$
We first claim that
$$\psi_1\psi_2\cdots\psi_n=\sum_{\tau\in U_{r,n}}\tau.$$
Clearly, $\psi_1=\sum_{\tau\in U_{r,1}}\tau$.
For the inductive step, suppose the claim is true for $n-1$ and identify the elements of $U_{r,n-1}$ with the set $\{\tau\in U_{r,n}:\,\tau_n=n\}$. 
Given such a permutation $\tau=\tau_1\cdots\tau_{n-1}\,n$, we have
\begin{equation}\label{eq:psi_algebra}
\tau\cdot\psi_n=\tau_1\cdots\tau_{n-1}\,n+n^{[1]}\tau_1\cdots\tau_{n-1}+\cdots+n^{[r-1]}\tau_1\cdots\tau_{n-1}.
\end{equation}
Moreover, it is obvious that $\tau\cdot\psi_n=\tau'\cdot\psi_n$ if and only if $\tau=\tau'$ and the induction is completed.

Define the linear mapping $\Psi:\mathbb{Z}[U_{r,n}]\to\mathbb{Z}[q,t]$ given by
$$\Psi(\tau) := t^{\sum_{z_i>0}(|\tau_i|-1)}q^{\col(\tau)}.$$
By construction, we have $\Psi(\psi_k)=1+t^{k-1}q+t^{k-1}q^2+\cdots+t^{k-1}q^{r-1}=1+t^{k-1}q[r-1]_q$ for $k\geq 1$.

We now prove that $\Psi(\psi_1\cdots\psi_n)=\Psi(\psi_1)\cdots\Psi(\psi_n)$.
Again, suppose for induction that the assertion is true for $n-1$.
It suffices to show that the following holds:
$$\Psi(\psi_1\cdots\psi_{n-1}\psi_{n})=\Psi(\psi_1\cdots\psi_{n-1})\Psi(\psi_n).$$
By \eqref{eq:psi_algebra} we have
\begin{align*}
\Psi(\tau\cdot\psi_n)
&=\Psi(\tau_1\cdots\tau_{n-1}n)+\Psi(n^{[1]}\tau_1\cdots\tau_{n-1})+\cdots+\Psi(n^{[r-1]}\tau_1\cdots\tau_{n-1})\\
&=\Psi(\tau)+t^{n-1}q\Psi(\tau)+\cdots+t^{n-1}q^{r-1}\Psi(\tau)\\
&=\Psi(\tau)\left(1+t^{n-1}q[r-1]_q\right)=\Psi(\tau)\Psi(\psi_n).
\end{align*}
Thus,
\begin{align*}
\Psi(\psi_1\cdots\psi_{n-1}\psi_{n})
&=\Psi\left(\sum_{\tau\in U_{r,n},\tau_n=n}\tau\cdot\psi_n\right)
=\sum_{\tau\in U_{r,n},\tau_n=n}\Psi(\tau\cdot\psi_n)\\
&=\Psi(\psi_n)\sum_{\tau\in U_{r,n},\tau_n=n}\Psi(\tau)
=\Psi(\psi_n)\Psi(\psi_1\cdots\psi_{n-1}),
\end{align*}
as desired.
Hence the result follows. \qed

%----length+lmaj------------------------------------------------
Now we are ready to present the first main result of this paper.

\begin{thm}\label{thm:G_length}
For $a=0,1$ and $b=0,1,\cdots, r-1$, we have the following signed Mahonian identities
\begin{equation}\label{eq:r_length}
\sum_{\pi\in G(r,n)} \chi_{a,b}(\pi)q^{\ell(\pi)} = \prod_{k=1}^n[k]_{(-1)^aq}\left(1+(-1)^{a(k-1)}\omega^bq^k[r-1]_{\omega^bq}\right),
\end{equation}
and
\begin{equation}\label{eq:r_lmaj}
\sum_{\pi\in G(r,n)} \chi_{a,b}(\pi)q^{\lmaj(\pi)} = \prod_{k=1}^n[k]_{(-1)^{a(k-1)}q} \left(1+(-1)^{a(k-1)}\omega^bq^k[r-1]_{\omega^bq}\right).
\end{equation}
\end{thm}
\proof
In what follows, let $(\Stat,\stat)=(\ell,\inv)$ or $(\lmaj,\maj)$.
Then,
\begin{eqnarray*}
&&\sum_{\pi\in G(r,n)} \chi_{a,b}(\pi)q^{\Stat(\pi)}\\
&=&\sum_{\pi\in G(r,n)}(-1)^{a(\inv(\pi)+\sum_{z_i>0}(|\pi_i|-1))}\omega^{b\,\col(\pi)}q^{\stat(\pi)+\sum_{z_i> 0}(|\pi_i|+z_i-1)}\\
&=&\sum_{\pi\in G(r,n)}\left((-1)^aq\right)^{\sum_{z_i>0}(|\pi_i|-1)} (\omega^bq)^{\col(\pi)}{(-1)^a}^{\inv(\pi)}q^{\stat(\pi)}\\
&=&\sum_{\tau\in U_{r,n}}\left((-1)^aq\right)^{\sum_{z_i>0}(|\tau_i|-1)}(\omega^bq)^{\col(\tau)}\sum_{\rho\in \mathfrak{S}_n}{(-1)^a}^{\inv(\rho)}q^{\stat(\rho)}\\
&=&F((-1)^aq,\omega^bq)\sum_{\rho\in \mathfrak{S}_n}{(-1)^a}^{\inv(\rho)}q^{\stat(\rho)}.
\end{eqnarray*}

Notice that if $\stat=\inv$, the last summation above turns out to be 
$\d\sum_{\rho\in \mathfrak{S}_n}\left((-1)^{a}q\right)^{\inv(\rho)}$. 
By \eqref{eq:A_Poincare}, the corresponding signed Mahonian identities for $\ell$ are
$$\prod_{k=1}^n[k]_{(-1)^aq}\left(1+(-1)^{a(k-1)}\omega^bq^k[r-1]_{\omega^bq}\right).$$
If $\stat=\maj$, that summation turns out to be the distribution of major index over $\mathfrak{S}_n$ when $a=0$ or the signed Mahonian for $\maj$ over $\mathfrak{S}_n$ when $a=1$. 
By \eqref{eq:A_Poincare} and \eqref{eq:A_signMaho}, the corresponding signed Mahonian polynomials for $\lmaj$ are
$$\prod_{k=1}^n[k]_{(-1)^{a(k-1)}q} \left(1+(-1)^{a(k-1)}\omega^bq^k[r-1]_{\omega^bq}\right).$$  \qed

%----link between lmaj and nmaj------------------------------------
%\begin{rmk}
%\begin{rm}
%It is easy to see that $\lmaj$ reduces to $\nmaj$ when $r=2$.
%Hence Theorem~\ref{thm:B_length}(1), (2) and (3) (resp. Theorem~\ref{thm:B_nmaj}(1), (2) and (3)) can be obtained by letting $r=2$ and $(a,b)=(1,1),(0,1),(1,0)$ in Equation~\eqref{eq:r_length} (resp. Equation~\eqref{eq:r_lmaj}) respectively.
%\end{rm}
%\end{rmk}

%----signed sorting-----------------------------------------------
\subsection{Signed Mahonian of the sorting index}\label{sec:Grn_sor}

Another way to represent elements of $G(r,n)$ is by means of ``bijections''.
Let $$\Sigma:=\{1,\ldots,n,\bar{1},\ldots,\bar{n},\bar{\bar{1}},\ldots,\bar{\bar{n}},\ldots,1^{[r-1]},\ldots,n^{[r-1]} \},$$ then $\pi=(z,\sigma)\in G(r,n)$ can be viewed as the bijection on $\Sigma$ such that $\pi(i)=\sigma_i^{[z_i]}$ for $i\in [n]$ and $\pi(\bar{i})=\overline{\pi(i)}$ for $i\in\Sigma$.
For example, $\bar{2}\bar{\bar{4}}1\bar{3}\bar{5}=(24135,(1,2,0,1,1))\in G(3,5)$ can be represented as the bijection
$$
\left(
\begin{array}{ccccccccccccccc}
\bar{\bar{1}} & \bar{\bar{2}} & \bar{\bar{3}} & \bar{\bar{4}} & \bar{\bar{5}} & \bar{1} & \bar{2} & \bar{3} & \bar{4} & \bar{5} & 1 & 2 & 3 & 4 & 5 \\
2 & \bar{4} & \bar{\bar{1}} & 3 & 5 & \bar{\bar{2}} & 4 & \bar{1} & \bar{\bar{3}} & \bar{\bar{5}} & \bar{2} & \bar{\bar{4}} & 1 & \bar{3} & \bar{5}
\end{array}
\right).
$$

For $1\leq i<j\leq n$ and $0\leq t<r$ let $(i^{[t]},j)$ be the transposition of swapping $\pi(i^{[t]})$ with $\pi(j)$, $\pi(i^{[t+1]})$ with $\pi(\bar{j})$, $\ldots$, $\pi(i^{[t+r-1]})$ with $\pi(j^{[r-1]})$. 
Also, for $1\leq i\leq n$ and $0<t<r$ let $(i^{[t]},i)$ be the action of adding $t$ bars on $\pi(i),\pi(\bar{i}),\ldots,\pi(i^{[r-1]})$.
For example, if $\pi=\bar{2}\bar{\bar{4}}1\bar{3}\bar{5}\in G(3,5)$, then $\pi\cdot(\bar{2},4)=\bar{2}314\bar{5}$ and $\pi\cdot(\bar{\bar{5}})=\bar{2}\bar{\bar{4}}1\bar{3}5$.
Formally, multiplying $\pi=\sigma_1^{[z_1]}\,\sigma_2^{[z_2]}\cdots\sigma_n^{[z_n]}\in G(r,n)$ on the right by $(i^{[t]},j)$, $i<j$, has the effect of replacing $\pi_j$ by $\sigma_i^{[z_i+t]}$ and $\pi_i$ by $\sigma_j^{[z_j-t]}$, while multiplying $\pi$ on the right by $(i^{[t]},i)$ has the effect of replacing $\pi_i$ by $\sigma_i^{[z_i+t]}$.

Given a permutation $\pi\in G(r,n)$, there is a unique expression $$\pi=(i_1^{[t_1]},j_1)(i_2^{[t_2]},j_2)\cdots(i_k^{[t_k]},j_k)$$ as a product of transpositions with $0<j_1<j_2<\cdots<j_k$ for some $k$.
The sorting index, $\sor(\pi)$, is given in \cite{ELW_15} by
\begin{equation}\label{eq:def_sor}
\sor(\pi) := \sum_{s=1}^k \left( j_s-i_s + \delta(t_s>0)\cdot\left( 2(i_s-1)+t_s \right) \right),
\end{equation}
where $\delta(\mathsf{A})=1$ if the statement $\mathsf{A}$ is true, or $\delta(\mathsf{A})=0$ otherwise.
For example, if $\pi=\bar{2}\bar{\bar{4}}1\bar{3}\bar{5}\in G(3,5)$, we have
\begin{align*}
\bar{2}\bar{\bar{4}}1\bar{3}\bar{5} \stackrel{~(\bar{\bar{5}},5)~}{\longrightarrow} 
\bar{2}\bar{\bar{4}}1\bar{3}5 \stackrel{~(\bar{2},4)~}{\longrightarrow}
\bar{2}3145 \stackrel{~(2,3)~}{\longrightarrow}
\bar{2}1345 \stackrel{~(\bar{\bar{1}},2)~}{\longrightarrow}
\bar{1}2345 \stackrel{~(\bar{\bar{1}},1)~}{\longrightarrow}
12345,
\end{align*}
and $\pi=(\bar{\bar{1}},1)(\bar{\bar{1}},2)(2,3)(\bar{2},4)(\bar{\bar{5}},5)$, so $\sor(\pi) = \ \big(1-1+(0+2)\big) + \big(2-1+(0+2)\big) + \big(3-2+0\big) + \big(4-2+(2\cdot 1+1)\big) + \big(5-5+(2\cdot 4+2)\big) = 21$.
%\begin{align*}
%\sor(\pi) = & \ \big(1-1+(0+2)\big) + \big(2-1+(0+2)\big) + \big(3-2+0\big) \\
%& + \big(4-2+(2\cdot 1+1)\big) + \big(5-5+(2\cdot 4+2)\big) = 21.
%\end{align*}
It has been proved in \cite{ELW_15} that $\sor$ has the same distribution with $\ell$ over $G(r,n)$ and thus is Mahonian, that is,
\begin{equation}\label{eq:equidist_sor}
\sum_{\pi\in G(r,n)}q^{\sor(\pi)}=\sum_{\pi\in G(r,n)}q^{\ell(\pi)}= \prod_{k=1}^n[k]_q(1+q^k[r-1]_q).
\end{equation}

To generate elements of $G(r,n)$, a simple way is to append the letter $n$ to the end of elements of $G(r,n-1)$, then pick an index $i$ and a color $t$ arbitrarily and apply the transposition $(i^{[t]},n)$.
Formally, we define the elements $\phi_i$ of the group algebra of $G(r,n)$ by $\phi_1:=1+(\bar{1},1)+\cdots+(1^{[r-1]},1)$ and for $2\leq k\leq n$ $$\phi_k:=1+\sum_{i=1}^{k-1}(i,k)+\sum_{t=1}^{r-1}\sum_{i=1}^{k}(i^{[t]},k).$$
Thus, for $\pi=\pi_1\cdots\pi_{n-1}n\in G(r,n)$ we have
\begin{align*}
\pi\cdot\phi_n = & \ \pi_1\cdots\pi_{n-1}n + \pi_1\cdots n\pi_{n-1} + \cdots + n\pi_2\cdots\pi_{n-1}\pi_1 \\
& + \pi_1\cdots\pi_{n-1}n^{[r-1]} + \pi_1\cdots n^{[r-1]}\pi_{n-1}^{[1]} + \cdots + n^{[r-1]}\pi_2\cdots\pi_{n-1}\pi_1^{[1]} \\
& + \cdots \\
& + \pi_1\cdots\pi_{n-1}n^{[1]} + \pi_1\cdots n^{[1]}\pi_{n-1}^{[r-1]} + \cdots + n^{[1]}\pi_2\cdots\pi_{n-1}\pi_1^{[r-1]}.
\end{align*}
%$$
%\begin{array}{rllll}
%\pi\cdot\phi_n & =\pi_1\cdots\pi_{n-1}n & +~\pi_1\cdots n\pi_{n-1} & +~\cdots & +~n\pi_2\cdots\pi_{n-1}\pi_1 \\
% & +~\pi_1\cdots\pi_{n-1}n^{[r-1]} & +~\pi_1\cdots n^{[r-1]}\pi_{n-1}^{[1]} & +~\cdots & +~n^{[r-1]}\pi_2\cdots\pi_{n-1}\pi_1^{[1]} \\
% & +~\cdots\cdots\cdots & & & \\ 
% & +~\pi_1\cdots\pi_{n-1}n^{[1]} & +~\pi_1\cdots n^{[1]}\pi_{n-1}^{[r-1]} & +~\cdots & +~n^{[1]}\pi_2\cdots\pi_{n-1}\pi_1^{[r-1]}.
%\end{array}
%$$
It has been proved in \cite{ELW_15} that 
\begin{equation}\label{eq:ELW_15}
\phi_1\phi_2\cdots\phi_n=\sum_{\pi\in G(r,n)}\pi.
\end{equation}

We are ready to derive the signed Mahonian polynomials with respect to the sorting index.

\begin{thm}\label{thm:sorting}
For $a=0,1$ and $b=0,1,\cdots, r-1$, we have the following signed Mahonian identities
%\begin{equation}\label{eq:r_sorting}
%\sum_{\pi\in G(r,n)} \chi_{a,b}(\pi)q^{\sor(\pi)} = \prod_{k=1}^n \left( 1 + (-1)^a[k-1]_q\left( q + \omega^{br}q^k[r-1]_q \right) + \omega^{b(r-t)}q^{2n-1}[r-1]_q \right).
%\end{equation}
\begin{align}
\sum_{\pi\in G(r,n)} & \chi_{a,b}(\pi)q^{\sor(\pi)} \notag \\ 
& = \prod_{k=1}^n \left( 1 + (-1)^a[k-1]_q\left( q + \omega^{br}q^k[r-1]_q \right) + \omega^b q^{2k-1}[r-1]_{\omega^b q} \right). \label{eq:r_sorting}
\end{align}
\end{thm}
\proof
%We proceed by induction.
For $a=0,1$ and $b=0,1,\cdots, r-1$, define the linear mapping $\Phi_{a,b}:\mathbb{Z}[G(r,n)]\to\mathbb{Z}[q]$ by $$\Phi_{a,b}(\pi):=\chi_{a,b}(\pi)q^{\sor(\pi)}.$$

Consider $\Phi_{a,b}(\phi_k)$ for $k\geq 1$.
By the definition of $\sor$, it is obvious that $\sor\big((i,k)\big)=k-i$ and $\sor\big((i^{[t]},k)\big)=k+i+t-2$ for $1\leq i\leq k$ and $1\leq t< r$.
Now, view $(i,k)=s_{k-1}s_{k-2}\cdots s_{i+1}s_is_{i+1}\cdots s_{k-2}s_{k-1}$.
Since $\chi_{a,b}$ is a 1-dim character, $$\chi_{a,b}\big((i,k)\big) = (-1)^{a((2k-2i-1)-0)}\cdot\omega^{b\cdot 0} = (-1)^a,$$
and thus 
\begin{equation}\label{eq:Phi_t=0}
\Phi_{a,b}\left( 1+\sum_{i=1}^{k-1}(i,k) \right) = 1 + \sum_{i=1}^{k-1}(-1)^aq^{k-i} = 1+(-1)^aq[k-1]_q.
\end{equation}
Similarly, since 
\begin{align*}
(i^{[t]},k)=\begin{cases}
s_{k-1}s_{k-2}\cdots s_1s_0^{r-t}s_1\cdots s_{i-2}s_{i-1}s_{i-2}\cdots s_1s_0^ts_1\cdots s_{k-2}s_{k-1}, & \text{if }1\leq i<k; \\
s_{k-1}s_{k-2}\cdots s_1s_0^{t}s_1\cdots s_{k-2}s_{k-1}, & \text{if }i=k,
\end{cases}
\end{align*}
%$$(i^{[t]},k)=s_{k-1}s_{k-2}\cdots s_1s_0^{r-t}s_1\cdots s_{i-2}s_{i-1}s_{i-2}\cdots s_1s_0^ts_1\cdots s_{k-2}s_k-1$$ for $1\leq i\leq k-1$ and $$(k^{[t]},k)=s_{k-1}s_{k-2}\cdots s_1s_0^{r-t}s_1\cdots s_{k-2}s_{k-1},$$ 
we have $$\chi_{a,b}\big((i^{[t]},k)\big) = (-1)^{a((2k+2i+r-5)-r)}\cdot\omega^{b\cdot r} = (-1)^a\omega^{br}, \quad \forall \ 1\leq i<k$$
and $$\chi_{a,b}\big((k^{[t]},k)\big) = (-1)^{a((2k+t-2)-t)}\cdot\omega^{b\cdot t} = \omega^{bt}.$$
Therefore, 
\begin{align}
\Phi_{a,b}\left( \sum_{t=1}^{r-1}\sum_{i=1}^{k}(i^{[t]},k) \right) & = \sum_{t=1}^{r-1}\left( \left( \sum_{i=1}^{k-1}(-1)^{a}\omega^{br}q^{k+i+t-2} \right) + \omega^{bt}q^{2k+t-2} \right) \notag \\
& = \sum_{t=1}^{r-1}\left( (-1)^{a}\omega^{br}q^{k+t-1}[k-1]_q + \omega^{bt}q^{2k+t-2} \right) \notag \\
& = (-1)^{a}\omega^{br}q^{k}[k-1]_q\sum_{t=1}^{r-1}q^{t-1} + q^{2k-1}\sum_{t=1}^{r-1}\omega^{bt}q^{t-1} \notag \\
& = (-1)^{a}\omega^{br}q^{k}[k-1]_q[r-1]_q + \omega^b q^{2k-1}[r-1]_{\omega^b q}. \label{eq:Phi_t>0}
\end{align}
Combining \eqref{eq:Phi_t=0} and \eqref{eq:Phi_t>0} yields
\begin{align*}
\Phi_{a,b}(\phi_k) = 1 + (-1)^a[k-1]_q\left( q + \omega^{br}q^k[r-1]_q \right) + \omega^b q^{2k-1}[r-1]_{\omega^b q}.
\end{align*}

By \eqref{eq:ELW_15} it suffices to show $\Phi_{a,b}(\phi_1\cdots\phi_n)=\Phi_{a,b}(\phi_1)\cdots\Phi_{a,b}(\phi_n)$.
Suppose for induction that the assertion is true for $n-1$.
We aim to show the following holds: $$\Phi_{a,b}(\phi_1\cdots\phi_{n-1}\phi_n)=\Phi_{a,b}(\phi_1\cdots\phi_{n-1})\Phi_{a,b}(\phi_n).$$
Let $\pi=\pi_1\cdots\pi_{n-1}n\in G(r,n)$.
Since $\pi(n)=n$, by the definition of $\sor$ we have $\sor(\pi\cdot(i,n))=\sor(\pi)+(n-i)$ for $1\leq i<n$, and $\sor(\pi\cdot(i^{[t]},n))=\sor(\pi)+(n+i+t-2)$ for $1\leq i\leq n$ and $1\leq t<r$.
Furthermore, by the fact that $\chi_{a,b}$ is an 1-dim character we have $\chi_{a,b}(\pi\cdot(i,n))=\chi_{a,b}(\pi)\cdot\chi_{a,b}((i,n))$ for $1\leq i<n$, and $\chi_{a,b}(\pi\cdot(i^{[t]},n))=\chi_{a,b}(\pi)\cdot\chi_{a,b}((i^{[t]},n))$ for $1\leq i\leq n$ and $1\leq t<r$.
Following a similar argument as in \eqref{eq:Phi_t=0} -- \eqref{eq:Phi_t>0} derives 
\begin{align*}
\Phi_{a,b}(\pi\cdot\phi_n) = \Phi_{a,b}(\pi)\Phi_{a,b}(\phi_n).
\end{align*}
Thus, 
\begin{align*}
\Phi_{a,b}(\phi_1\cdots\phi_{n-1}\phi_n) & = \Phi_{a,b}\left( \sum_{\pi\in G(r,n),\pi(n)=n}\pi\cdot\phi_n \right) = \sum_{\pi\in G(r,n),\pi(n)=n} \Phi_{a,b}(\pi\cdot\phi_n) \\
& = \Phi_{a,b}(\phi_n) \sum_{\pi\in G(r,n),\pi(n)=n} \Phi_{a,b}(\pi) = \Phi_{a,b}(\phi_n)\Phi_{a,b}(\phi_1\cdots\phi_{n-1}),
\end{align*}
as desired. 
\qed

%%%%%%%%%%%%%%%%%%%%%%%%%%%%%%%%%%%%%%%%%%%%%%%%%%%%%%%%%%%%%%%%%%%%%%%%%%%
%%%%%%%%%%%%%%%%%%%%%%%%%%%%%%%%%%%%%%%%%%%%%%%%%%%%%%%%%%%%%%%%%%%%%%%%%%%
\section{Signed Mahonian on $B_n$}\label{sec:Bn}
For Coxeter groups of type $B_n=G(2,n)$ there are some Mahonian statistics other than $\ell_B$ and $\lmaj$.
We take a close look in this section.

\subsection{Signed permutation group}\label{Sec:Bn_Bn}

Let $B_n$ be the \emph{signed permutation group} of $\{1,2,\ldots,n\}$, which consists of all bijections $\pi$ of $\{\pm 1,\ldots,\pm n\}$ onto itself such that $\pi(-i)=-\pi(i)$. 
%For simplicity we denote $\pi_i:=\pi(i)$. 
Elements in $B_n$ are centrally symmetric and hence can be simply written in the form $\pi=\pi_1\pi_2\cdots\pi_n$, where $\pi_i:=\pi(i)$.
%Throughout this paper we will only discuss window notations. 
$B_n$ is also the Coxeter group of type $B_n$ with generators $s_0=(1,-1)$ and $s_i=(i,i+1)$ for $i=1,\ldots,n-1$. 
Let $\ell_B$ be the length function of $B_n$ with respect to this set of generators. 
Recall that 
$$\ell_B(\pi)=\inv(\pi)-\sum_{i\in\Nega(\pi)}\pi_i.$$
By denoting $-i$ by $i^{[1]}$, it is obvious that $B_n=G(2,n)$, and $\ell_B$ turns out to be a special case of $\ell$ with $r=2$ in \eqref{eq:length_Grn}.

The search of a Mahonian major statistic on $B_n$ was a long-standing problem and was first solved by Adin and Roichman~\cite{Adin_Roichman_01} by introducing the flag major index, $\fmaj$.
So far there are at least three type $B_n$ Mahonian statistics in the literature, defined as follows.
\begin{itemize}
\item $\fmaj$. The \emph{flag major index}~\cite{Adin_Roichman_01} is defined by
$$\fmaj(\pi):=2\cdot\maj_F(\pi)+\nega(\pi).$$
\item $\Fmaj$. The \emph{F-major index}~\cite{Adin_Roichman_01} is defined by
$$\Fmaj(\pi):=2\cdot\maj(\pi)+\nega(\pi).$$
\item $\nmaj$. The \emph{negative major index}~\cite{Adin_Brenti_Roichman_01} is defined by
$$\nmaj(\pi):=\maj(\pi)-\sum_{i\in \Nega(\pi)}\pi_i.$$
\end{itemize}
Here $\maj$ is defined as before, while $\maj_F$ is computed similarly but with respect to the \emph{flag-order}:
$$\bar{1}<\cdots<\bar{n}<1<\cdots<n.$$
For example, if $\pi=\bar{3}1\bar{6}2\bar{5}\bar{4}$, then $\fmaj(\pi)=2\cdot 11+4=26$, $\Fmaj(\pi)=2\cdot 6+4=16$, and $\nmaj(\pi)=6-(-18)=24$.
Note that $\nmaj$ is exactly $\lmaj$ with $r=2$.

In addition to major-type statistics on $B_n$, the type $B_n$ sorting index $\sor_B$ is defined by Petersen~\cite{Petersen_11}.
It's definition can be reduced from \ref{eq:def_sor} by simply plugging $r=2$.

To summarize, we from \eqref{eq:equidist_Grn} have
\begin{equation}\label{eq:B_equidist}
\sum_{\pi\in B_n}q^{\ell_B(\pi)}=\sum_{\pi\in B_n}q^{\fmaj(\pi)}
=\sum_{\pi\in B_n}q^{\Fmaj(\pi)}=\sum_{\pi\in B_n}q^{\nmaj(\pi)}=\sum_{\pi\in B_n}q^{\sor_B(\pi)}=\prod_{k=1}^{n}[2k]_q.
\end{equation}

%%%%%%%%%%%%%%%%%%%%%%%
\subsection{Signed Mahonian on $B_n$}\label{Sec:Bn_SM}
The four 1-dim characters of $B_n$, say the trivial character 1, $(-1)^{\nega(\pi)}$, $(-1)^{\ell_B(\pi)}$ and $(-1)^{\inv(|\pi|)}$, can be obtained from Lemma~\ref{lem:1-dim} by letting $r=2$.
More precisely,
\begin{align*}
\chi_{0,0}(\pi) &= 1, \\
\chi_{0,1}(\pi) &= (-1)^{\nega(\pi)}, \\
\chi_{1,1}(\pi) &= (-1)^{\ell_B(\pi)}, \\
\chi_{1,0}(\pi) &= (-1)^{\ell_B(\pi)-\nega(\pi)} = (-1)^{\inv(\pi)-\big(\sum_{i\in\Nega(\pi)}\pi_i\big) - \nega(\pi)} \\
&= (-1)^{\inv(\pi)+\sum_{i\in\Nega(\pi)}(|\pi_i|-1)} \stackrel{(*)}{=} (-1)^{\inv(|\pi|)},
\end{align*}
where the proof of (*) is relegated to Lemma~\ref{lem:inv}.
%and $(a,b)=(0,0), (0,1), (1,0), (1,1)$, which are respectively denoted herein by the trivial character 1, $(-1)^{\nega(\pi)}$, $(-1)^{\inv(|\pi|)}$ and $(-1)^{\ell_B(\pi)}$. 
The signed Mahonian polynomials of $\fmaj$ together with the four 1-dim characters were obtained by Adin et al.~\cite{Adin_Gessel_Roichman_05}, and the signed Mahonian polynomials of $\ell_B$ or $\nmaj$ can be derived directly from Theorem~\ref{thm:G_length} as follows.

\begin{cor}
For $\ell_B$ we have the following signed Mahonian polynomials:
\begin{enumerate}
%\item $\d\sum_{\pi\in B_n}1\,q^{\ell^B(\pi)}=\prod_{k=1}^n[2k]_{q}.$
\item $\d\sum_{\pi\in B_n}(-1)^{\ell_B(\pi)}q^{\ell_B(\pi)}=\prod_{k=1}^n[2k]_{-q}$;
\item $\d\sum_{\pi\in B_n}(-1)^{\nega(\pi)}q^{\ell_B(\pi)}=\prod_{k=1}^n[2]_{-q^k}[k]_q$;
\item $\d\sum_{\pi\in B_n}(-1)^{\inv(|\pi|)}q^{\ell_B(\pi)}=\prod_{k=1}^n[2]_{(-1)^{k-1}q^k}[k]_{-q}$.
%\item $\d\sum_{\pi\in B_n}(-1)^{\inv_A(\pi)}q^{\ell^B(\pi)}=\prod_{k=1}^n[2]_{q^k }[k]_{-q}$.
\end{enumerate}
For $\nmaj$ we have the following signed Mahonian polynomials:
\begin{enumerate}
%\item $\d\sum_{\pi\in B_n}1\,q^{\nmaj(\pi)}=\prod_{k=1}^n[2k]_{q}.$
\item $\d\sum_{\pi\in B_n}(-1)^{\ell_B(\pi)}q^{\nmaj(\pi)} =\prod_{k=1}^n[2]_{(-q)^k}[k]_{(-1)^{k-1}q}$;
\item $\d\sum_{\pi\in B_n}(-1)^{\nega(\pi)}q^{\nmaj(\pi)}=\prod_{k=1}^n[2]_{-q^k}[k]_q$;
\item $\d\sum_{\pi\in B_n}(-1)^{\inv(|\pi|)}q^{\nmaj(\pi)}=\prod_{k=1}^n[2]_{(-1)^{k-1}q^k}[k]_{(-1)^{k-1}q}$.
%\item $\d\sum_{\pi\in B_n}(-1)^{\inv_A(\pi)}q^{\nmaj(\pi)}=\prod_{k=1}^n[2]_{q^k}[k]_{(-1)^{k-1}q}$.
\end{enumerate}
\end{cor}

%It must be noted that the cases $\Fmaj$ and $\nmaj$ for $\chi(\pi)=\sign(\pi)$ were derived in \cite[Theorem 5.7, 5.8]{Fire_08}.

The signed Mahonian polynomial of $\sor_B$ can be obtained from Theorem~\ref{thm:sorting} as follows.
\begin{cor}
For $\sor_B$ we have the following signed Mahonian polynomials:
\begin{enumerate}
\item $\d\sum_{\pi\in B_n}(-1)^{\ell_B(\pi)}q^{\sor_B(\pi)}=\prod_{k=1}^n \left( 1-q[2k-1]_q \right)$;
\item $\d\sum_{\pi\in B_n}(-1)^{\nega(\pi)}q^{\sor_B(\pi)}=\prod_{k=1}^n \left( [2k-1]_q-q^{2k-1} \right)$;
\item $\d\sum_{\pi\in B_n}(-1)^{\inv(|\pi|)}q^{\sor_B(\pi)}=\prod_{k=1}^n \left( 1-q[2k-2]_q+q^{2k-1} \right)$.
\end{enumerate}
\end{cor}

In the rest of this subsection we focus on the signed Mahonian polynomials of $\Fmaj$. 
Following the same fashion as in $G(r,n)$, we consider the decomposition $B_n=U_n\cdot\mathfrak{S}_n$, where $U_n=U_{2,n}=\{\tau\in B_n| \tau_1<\tau_2<\cdots<\tau_n\}$.
It is obvious that $\inv(\tau\rho)=\inv(\rho)$, $\maj(\tau\rho)=\maj(\rho)$, and $\Nega(\tau\rho)=\Nega(\tau)$ for $\tau\in U_n$ and $\rho\in\mathfrak{S}_n$.
By letting $r=2$ in Lemma~\ref{lem:r_bivariate} we obtain that
\begin{equation}\label{eq:2_F(t,q)}
F(t,q)\Big|_{r=2} = \sum_{\tau\in U_n}t^{\sum_{i\in\Nega(\tau)}(|\tau_i|-1)}q^{\nega(\tau)} = [2]_{t^{k-1}q}.
\end{equation}
Furthermore, by a similar argument as in Lemma~\ref{lem:r_bivariate}, we have
\begin{equation}\label{eq:B(t,q)}
B(t,q):=\sum_{\tau\in U_n}t^{\nega(\tau)}q^{-\sum_{i\in\Nega(\tau)}\tau_i}=\prod_{k=1}^n[2]_{tq^k}.
\end{equation}

%We decompose $B_n$ into $U_n\cdot \mathfrak{S}_n$, where $U_n:=\{\tau\in B_n| \tau_1<\tau_2<\cdots<\tau_n\}$. 
%It is clear that
%$$B_n=\biguplus_{\tau\in U_n}\{\tau\rho:\,\rho\in \mathfrak{S}_n\}=\biguplus_{\rho\in \mathfrak{S}_n}\{\tau\rho:\,\tau\in U_n\}.$$
%We omit the proofs of the following simple facts.
%\begin{lem}\label{lem:B_decom}
%For any $\tau\in U_n$ and $\rho\in \mathfrak{S}_n$, we have $\inv(\tau\rho)=\inv(\rho)$, 
%$\maj(\tau\rho)=\maj(\rho)$, $\nega(\tau\rho)=\nega(\tau)$, and
%$\sum_{i\in\Nega(\tau\rho)}(\tau\rho)_i=\sum_{i\in\Nega(\tau)}\tau_i.$
%\end{lem}

The following lemma is essential to our main results in this section.

%\begin{lem}\label{lem:n&s} We have
%\begin{equation}\label{eq:B(t,q)}
%B(t,q):=\sum_{\tau\in U_n}t^{\nega(\tau)}q^{\ell_B(\tau)}=\prod_{k=1}^n[2]_{tq^k},
%\end{equation}
%and
%\begin{equation}\label{eq:B_evenNeg}
%\sum_{\tau\in U_n}(-1)^{|\Nega(\tau)\cap\{i:\,\tau_i\text{ is even}\}|}q^{\nega(\tau)} = \prod_{k=1}^n[2]_{(-1)^{k-1}q}.
%\end{equation}
%\end{lem}
%\proof 
%We first consider \eqref{eq:B(t,q)}. 
%For $k\ge 1$ let $\psi_k=1+s_{k-1}s_{k-2}\cdots s_0$. 
%By the same argument in the proof of Lemma~\ref{lem:r_bivariate} we have
%$$\psi_1\psi_2\cdots\psi_n=\sum_{\tau\in U_n}\tau.$$
%
%Now, define the linear mapping $\Psi:\mathbb{Z}[U_n]\to\mathbb{Z}[q,t]$ given by
%$$\Psi(\tau) = t^{\nega(\tau)}q^{\ell_B(\tau)}.$$
%By construction, we have $\Psi(\psi_k)=1+tq^k=[2]_{tq^k}$. 
%It is easy to see that $\Psi(\psi_1\cdots\psi_n)=\Psi(\psi_1)\cdots\Psi(\psi_n)$ and hence
%\eqref{eq:B(t,q)} follows. 
%
%Similarly, if we define the linear mapping $\Theta:\mathbb{Z}[U_n]\to\mathbb{Z}[q,t]$ by 
%$$\Theta(\tau) = (-1)^{|\Nega(\tau)\cap\{i:\,\tau_i\text{ is even}\}|}q^{\nega(\tau)},$$
%then \eqref{eq:B_evenNeg} can be derived by the same construction. \qed

\begin{lem}\label{lem:inv}
For any signed permutation $\pi\in B_n$, we have
$$(-1)^{\inv(|\pi|)}=(-1)^{\inv(\pi)+\sum_{i\in\Nega(\pi)}(|\pi_i|-1)}.$$
\end{lem}
\proof 
If $\Nega(\pi)=\emptyset$, then $|\pi|=\pi$ and the identity trivially holds.
So we only need to consider the case of $\nega(\pi)>0$.
Assume $\pi_k=\max\{|\pi_i| : i\in\Nega(\pi)\}$, and let $\pi'=\pi_1\cdots\pi_{k-1}|\pi_k|\pi_{k+1}\cdots\pi_n$.
For example, if $\pi=\bar{2}73\bar{1}6\bar{5}4$, then $\pi_k=\bar{5}$ and $\pi'=\bar{2}73\bar{1}654$.

Let $s=|\{\pi_i: i<k \text{ and }|\pi_i|<|\pi_k|\}|$ and $t=|\{\pi_i: i>k \text{ and } |\pi_i|<|\pi_k|\}|$. One has $s+t=|\pi_k|-1$.
By a simple calculation, we have
$$\inv(\pi')=\inv(\pi)-s+t=\inv(\pi)+|\pi_k|-1-2s,$$ 
implying $(-1)^{\inv(\pi')}=(-1)^{\inv(\pi)+|\pi_k|-1}.$
%Hence, by the same arguments we have
%\begin{align}
%(-1)^{\inv(|\pi|)+\inv(\pi)} &= (-1)^{\sum_{i\in\Nega(\pi)}(|\pi_i|-1)} \label{eq:proof_Bn_1} \\
%&= (-1)^{|\Nega(\pi)\cap\{i: \pi_i\text{ is even}\}|}.\notag
%\end{align}
Hence, the result follows by considering iteratively all elements in $\Nega(\pi')$.
\qed

We are ready for the main result in this section.
Note that Theorem~\ref{thm:B_Fmaj}(1) has appeared in~\cite[Theorem 5.7]{Fire_08}.
We still put it here for the sake of completeness.
%----Fmaj----------------------------------------------------

\begin{thm}\label{thm:B_Fmaj}
For $\Fmaj$ we have the following signed Mahonian polynomials.
\begin{enumerate}
\item $\d\sum_{\pi\in B_n}(-1)^{\ell_B(\pi)}q^{\Fmaj(\pi)}=\prod_{k=1}^n[2]_{(-1)^kq}[k]_{(-1)^{k-1}q^2}$.
\item $\d\sum_{\pi\in B_n}(-1)^{\nega(\pi)}q^{\Fmaj(\pi)} =\prod_{k=1}^n[2]_{-q}[k]_{q^2}$.
\item $\d\sum_{\pi\in B_n}(-1)^{\inv(|\pi|)}q^{\Fmaj(\pi)} =\prod_{k=1}^n[2]_{(-1)^{k-1}q}[k]_{(-1)^{k-1}q^2}$.
%\item $\d\sum_{\pi\in B_n}(-1)^{\inv_A(\pi)}q^{\Fmaj(\pi)}=\prod_{k=1}^n[2]_q[k]_{(-1)^{k-1}q^2}$.
\end{enumerate}
\end{thm}
\proof 
Recall that $\Fmaj(\pi)=2\cdot\maj(\pi)+\nega(\pi)$ and $B_n=U_n\cdot\mathfrak{S}_n$.

(1) By \eqref{eq:A_signMaho} and \eqref{eq:B(t,q)},
the left-hand side is equal to
\begin{eqnarray*}
&&\sum_{\tau\in U_n}(-1)^{-\sum_{i\in\Nega(\tau)}\tau_i}q^{\nega(\tau)} \sum_{\rho\in \mathfrak{S}_n}(-1)^{\inv(\rho)}q^{2\cdot\maj(\rho)}\\
&=& B(q,-1)\prod_{k=1}^n[k]_{(-1)^{k-1}q^2} = \prod_{k=1}^n[2]_{(-1)^kq}[k]_{(-1)^{k-1}q^2}.
\end{eqnarray*}

(2) By \eqref{eq:A_Poincare} and \eqref{eq:B(t,q)}, the left-hand side is equal to
$$\sum_{\tau\in U_n}(-q)^{\nega(\tau)} \sum_{\rho\in \mathfrak{S}_n}q^{2\cdot\maj(\rho)}
= B(-q,1)\prod_{k=1}^n[k]_{q^2} = \prod_{k=1}^n[2]_{-q}[k]_{q^2}.$$

(3) By \eqref{eq:A_signMaho}, \eqref{eq:2_F(t,q)} and Lemmas~\ref{lem:inv}, the left-hand side is equal to
\begin{eqnarray*}
&&\sum_{\pi\in B_n}(-1)^{\inv(\pi)+\sum_{i\in\Nega(\pi)}(|\pi_i|-1)}q^{2\cdot\maj(\pi)+\nega(\pi)}\\
&=& \sum_{\tau\in U_n}(-1)^{\sum_{i\in\Nega(\tau)}(|\tau_i|-1)}q^{\nega(\tau)} \sum_{\rho\in \mathfrak{S}_n}(-1)^{\inv(\rho)}q^{2\cdot\maj(\rho)}\\
&=& F(-1,q)\Big|_{r=2}\prod_{k=1}^n[k]_{(-1)^{k-1}q^2} = \prod_{k=1}^n[2]_{(-1)^{k-1}q}[k]_{(-1)^{k-1}q^2}.
\end{eqnarray*}
\qed

\section{Signed Mahonian on $D_n$}\label{sec:Dn}

\subsection{Even-signed permutation group}\label{Sec:Dn_Dn}

The \emph{even-signed permutation group} $D_n$ is the subgroup of $B_n$ defined by 
$$D_n:=\{\pi\in B_n:\, \nega(\pi)\text{ is even}\},$$ 
which consists of those permutations with an even number of negatives among $\pi_1,\ldots,\pi_n$. 
%Note that $D_n=G(2,2,n)$.
$D_n$ is known as the ``Coxeter group of type $D_n$'' which has the generators $s_0',s_1,\ldots,s_{n-1}$, where $s_0'=(\bar{1},2)$ and $s_i=(i,i+1)$ for $i\geq 1$. 
Let $\ell_D$ be the corresponding length function of $D_n$. 
It is known~\cite{Bjorner_Brenti_05, Humphreys_90} that the generating function for the distribution of $\ell_D$ is
\begin{equation}\label{eq:D_equidist}
\sum_{\pi\in D_n} q^{\ell_D(\pi)}=[n]_q\prod_{k=1}^{n-1}[2k]_q=[2]_q[4]_q\cdots [2n-2]_q[n]_q
\end{equation}
and a combinatorial description of $\ell_D$ is
$$ \ell_D(\pi)=\inv(\pi)-\nega(\pi)-\sum_{i\in\Nega(\pi)}\pi_i.$$

Biagioli \cite{Biagioli_03} defined the index
$$\dmaj(\pi):=\maj(\pi)-\nega(\pi)-\sum_{i\in \Nega(\pi)}\pi_i$$
on $D_n$ and proved that it is Mahonian.

$D_n$ can also be generated by $$\mathcal{T}_n:=\{t_{i,j}:\,1\leq|i|<j\leq n\} \cup \{t_{\bar{i},i}:\,1<i\leq n\},$$ where $t_{i,j}=(i,j)$ for $1\leq|i|<j\leq n$ and $t_{\bar{i},i}=(\bar{1},1)(\bar{i},i)$ for $1<i\leq n$.
Then, any $\pi\in D_n$ has a unique factorization in the form $$\pi=t_{i_1,j_1}t_{i_2,j_2}\cdots t_{i_k,j_k}$$ with $0<j_1<j_2<\cdots<j_k\leq n$ for some $k$.
Petersen \cite{Petersen_11} defined the sorting index of type $D_n$ by
\begin{align*}
\sor_D(\pi):=\sum_{s=1}^{k} \left( j_s-i_s - 2\delta(i_s<0) \right)
\end{align*}
and proved that it is Mahonian.
For example,
$\pi=\bar{3}24\bar{5}1=t_{\bar{1},3}t_{3,4}t_{\bar{4},5}$ has the sorting index $\sor_D(\pi)= (3-(-1)-2)+ (4-3)+ (5-(-4)-2)= 10$.

%%%%%%%%%%%%%%%%%%%%%%%%%%%
\subsection{Signed Mahonian on $D_n$}\label{sec:Dn_SM}
The group $D_n$ has two 1-dim characters~\cite{Reiner_95}: the trivial character $1$ and $(-1)^{\ell_D(\pi)}$. 
In this subsection we derive the signed Mahonian polynomials $\sum_{\pi\in D_n } (-1)^{\ell_D(\pi)} q^{\stat(\pi)}$, where $\stat=\ell_D,\dmaj$ or $\sor_D$.
%It must be noted that we do not take $(-1)^{\nega(\pi)}$ and $\sign(|\pi|)$ into consideration since the former is identically 1 and the latter is equal to $\sign(\pi)$ for $\pi\in D_n$.

%----preparation------------------------------------------------
\medskip

We first consider $\stat=\ell_D$ or $\dmaj$.
Similar to the case of $B_n$, we consider the decomposition $D_n=U^D_n\cdot \mathfrak{S}_n$, where $U^D_n:=\{\pi\in U_n:\,\nega(\pi)\text{ is even}\}$. 
%Let $U^D_n:=\{\pi\in U_n:\,\nega(\pi)\text{ is even}\}$. 
%Then one has $D_n=U^D_n\cdot \mathfrak{S}_n$ and
%\begin{equation}\label{eq:D_decom}
%D_n = \biguplus_{\tau\in U_n^D}\{\tau\rho:\, \rho\in \mathfrak{S}_n\}=\biguplus_{\rho\in \mathfrak{S}_n}\{\tau\rho:\, \tau\in U_n^D\}.
%\end{equation}

\begin{lem}\label{lem:D_stat}
We have
$$B^D(q):=\sum_{\tau\in U^D_n}q^{\ell_D(\tau)}=\prod_{k=1}^{n-1}[2]_{q^k}.$$
%$$ B^D(q):=\sum_{\substack{T\subseteq [n],\\|T|\text{ is even}}} q^{\left(\sum_{i\in T}i\right)-|T|}=\prod_{k=1}^{n-1}[2]_{q^k}.$$
\end{lem}
\proof 
Let $[n]:=\{1,2,\ldots,n\}$. 
Then $\sum_{\tau\in U^D_n}q^{\ell_D(\tau)}$ is equal to
\begin{eqnarray*}
&&\sum_{\substack{T\subseteq [n],\\|T|\text{ is even}}} q^{\left(\sum_{i\in T}i\right)-|T|}
=\sum_{\substack{T\subseteq \{0,1,\cdots,n-1\},\\|T|\text{ is even}}} q^{\sum_{i\in T}i}\\
&=&\sum_{\substack{T\subseteq [n-1],\\|T|\text{ is odd}}}q^{\sum_{i\in T\cup\{0\}}i}
+\sum_{\substack{T\subseteq [n-1],\\|T|\text{ is even}}}q^{\sum_{i\in T}i}
=\sum_{T\subseteq [n-1]}q^{\sum_{i\in T}i}
=\prod_{k=1}^{n-1}(1+q^k),
\end{eqnarray*}
as desired. \qed

%----length----------------------------------------------

%We are ready to state the objective signed Mahonian polynomials.

\begin{thm}\label{D_length&dmaj}
For $D_n$ we have the following signed Mahonian polynomials.
\begin{enumerate}
\item $\d\sum_{\pi\in D_n}(-1)^{\ell_D(\pi)}q^{\ell_D(\pi)}=[n]_{-q}\prod_{k=1}^{n-1}[2k]_{-q}.$
%\item $\d\sum_{\pi\in D_n}(-1)^{\nega(\pi)}q^{\ell^D(\pi)}=[n]_{q}\prod_{k=1}^{n-1}[2k]_{q}.$
%\item $\d\sum_{\pi\in D_n}\sign(|\pi|)q^{\ell^D(\pi)}=[n]_{-q}\prod_{k=1}^{n-1}[2]_{(-q)^k}[k]_{-q}$.
%\item $\d\sum_{\pi\in D_n}(-1)^{\inv_A(\pi)}q^{\ell^D(\pi)}=[n]_{-q}\prod_{k=1}^n[2]_{q^k}[k]_{-q}$.
\item $\d\sum_{\pi\in D_n}(-1)^{\ell_D(\pi)}q^{\dmaj(\pi)}=[n]_{(-1)^{n-1}q}\prod_{k=1}^{n-1}[2]_{(-q)^k}[k]_{(-1)^{k-1}q}.$
%\item $\d\sum_{\pi\in D_n} (-1)^{\inv_A(\pi)}q^{\dmaj(\pi)}=[n]_{(-1)^{n-1}q}\prod_{k=1}^{n-1}[2]_{q^k}[k]_{(-1)^{k-1}q}.$
\end{enumerate}
\end{thm}
\proof
(1) is obtained by replacing $q$ by $-q$ in \eqref{eq:D_equidist}.

(2) Recall that $D_n=U^D_n\cdot \mathfrak{S}_n$. 
By Lemma~\ref{lem:D_stat} and \eqref{eq:A_signMaho}, the left-hand side is equal to
$$\sum_{\tau\in U^D_n}(-q)^{-\left(\sum_{i\in\Nega(\tau)}\tau_i\right)-\nega(\tau)} \sum_{\rho\in \mathfrak{S}_n}(-1)^{\inv(\rho)}q^{\maj(\rho)}=B^D(-q)\prod_{k=1}^n [k]_{(-1)^{k-1}q},$$
as desired. \qed
%(2) By Lemma~\ref{lem:B_decom} and \eqref{eq:D_decom}, the left-hand side is equal to
%$$\sum_{\tau\in U^D_n}q^{-\left(\sum_{i\in\Nega(\rho)}\tau_i\right)-\nega(\tau)}\sum_{\rho\in \mathfrak{S}_n}(-q)^{\inv(\rho)} = B^D(q)\prod_{k=1}^{n}[k]_{-q} = [n]_{-q}\prod_{k=1}^n[2]_{q^k}[k]_{-q}.$$ 
%\qed

\noindent \emph{Remark.}
Note that Biagioli and Caselli \cite{Biagioli_Caselli_04} defined two other Mahonian indices over $D_n$, namely the flag major index $\fmaj$ and $D$-major index $\Dmaj$. 
The signed Mahonian polynomials of $\fmaj$ and $\Dmaj$ were derived in \cite{Biagioli_Caselli_12} and \cite{Biagioli_06}, respectively.
%as $$\sum_{\pi\in D_n } \sign(\pi)q^{\Dmaj(\pi)}=[n]_q\prod_{k=1}^{n-1}[2k]_{(-1)^{k}q},$$ which does not overlap our results.

In what follows we consider $\sor_D$.
Define elements $\theta_i$ of the group algebra of $D_n$ by $\theta_1:=1$ and for $2\leq k\leq n$ $$\theta_k:=1+\sum_{i=1}^{k-1}t_{i,k}+\sum_{i=1}^{k}t_{\bar{i},k}.$$
It has been shown in \cite{Petersen_11} that 
\begin{equation}\label{eq:D_algebra}
\theta_1\theta_2\cdots\theta_n = \sum_{\pi\in D_n}\pi.
\end{equation}

\begin{thm}\label{D_sor}
For $\sor_D$ we have the following signed Mahonian polynomials:
$$\sum_{\pi\in D_n}(-1)^{\ell_D(\pi)}q^{\sor_D(\pi)} = \prod_{k=2}^n \left( 1-(q+q^{k-1})[k-1]_q+q^{2k-2} \right).$$
\end{thm}
\proof
Define the linear mapping $\Theta:\mathbb{Z}[D_n]\to\mathbb{Z}[q]$ by $$\Theta(\pi):=(-1)^{\ell_D(\pi)}q^{\sor_D(\pi)}.$$
We first claim that, for $k\geq 2$,
\begin{equation}\label{eq:Theta}
\Theta(\theta_k) = 1-(q+q^{k-1})[k-1]_q+q^{2k-2}.
\end{equation}
Note that $\Theta(\theta_1)=1$.

By the definition of $\sor_D$, one has $\sor_D(t_{i,k})=k-i$ for $1\leq i<k$, and $\sor_D(t_{\bar{i},k})=k+i-2$.
Following the same argument as in \eqref{eq:Phi_t=0}, we have 
\begin{align*}
\Theta\left( 1+\sum_{i=1}^{k-1}t_{i,k} \right) = 1-q[k-1]_q.
\end{align*}
Next, observe that
\begin{align*}
t_{\bar{i},k} = \begin{cases}
s_{i-1}s_{i-2}\cdots s_1s_{k-1}s_{k-2}\cdots s_2 s'_0 s_2 \cdots s_{k-2} s_{k-1} s_1\cdots s_{i-2} s_{i-1}, & \text{if } 1\leq i<k; \\
s_{k-1}s_{k-2}\cdots s_2 s'_0 s_1 s_2 \cdots s_{k-2} s_{k-1}, & \text{if } i=k.
\end{cases}
\end{align*}
Since generators $s_0',s_1,\ldots,s_{n-1}$ are involutions, the expression of $t_{\bar{i},k}$ above yields
\begin{align*}
(-1)^{\ell_D(t_{\bar{i},k})} = \begin{cases}
(-1)^{2(i-1+k-2)+1}=-1, & \text{if } 1\leq i<k; \\
(-1)^{2(k-1)}=1, & \text{if } i=k.
\end{cases}
\end{align*}
This implies that $\Theta(t_{\bar{i},k})=-q^{k+i-2}$ for $1\leq i<k$, and $\Theta(t_{\bar{k},k})=q^{2k-2}$.
%by the fact that $(-1)^{\ell_D(\pi)}$ is an 1-dim character on $D_n$ we have 
%$$\Theta(t_{\bar{i},k})=(-1)^{2(i-1+k-2)+1}q^{k+i-2}=-q^{k+i-2}, \quad \forall \ 1\leq i<k$$ and $$\Theta(t_{\bar{k},k})=(-1)^{2(k-1)}q^{2k-2}=q^{2k-2}.$$
Therefore, 
\begin{align*}
\Theta(\theta_k) & = 1 - q[k-1]_q - \sum_{i=1}^{k-1}q^{k+i-2} + q^{2k-2}, %\\
%& = 1 - q[k-1]_q - q^{k-1}[k-1]_q + q^{2k-2},
\end{align*}
and thus \eqref{eq:Theta} follows.

By \eqref{eq:D_algebra} it suffices to show $\Theta(\theta_1\cdots\theta_n)=\Theta(\theta_1)\cdots\Theta(\theta_n)$.
We omit the rest of the proof since it can be dealt with in the same fashion as we did in the proof of Theorem~\ref{thm:sorting}.
\qed

%%%%%%%%%%%%%%%%%%%%%%%%%%%%%%%%%%%%%%%%%%%%%%%%%%%%%%%%%%%%%%%%%%%
%%%%%%%%%%%%%%%%%%%%%%%%%%%%%%%%%%%%%%%%%%%%%%%%%%%%%%%%%%%%%%%%%%%

\section{Signed counting on another class of statistics on $G(r,n)$}\label{sec:others}

In this section we investigate four equidistributed statistics on $G(r,n)$: $\fmaj$, $\rmaj$, $\fmaf$ and $\rinv$, the last of which is new. 
Hence we can consider the polynomial
$$\sum_{\pi\in G(r,n)}(-1)^{\stat_1(\pi)}q^{\stat_2(\pi)}$$
by taking two of these four statistics.
It turns out that among the possible twelve cases in many of which we can have nice closed forms.

%%%%%%%%%%%%%%%%%%%%%%%
\subsection{Equidistributed statistics}\label{sec:others_intro}
Throughout this section both $\inv$ and $\maj$ are calculated with respect to the linear order
$$(1^{[r-1]}<\cdots<n^{[r-1]})<\cdots<(1^{[1]}<\cdots<n^{[1]})<(1<\cdots<n).$$

Let $\pi=(z, \sigma)\in G(r,n)$ and define the following statistics:
\begin{itemize}
\item $\fmaj$. The \emph{flag major index}~\cite{Adin_Roichman_01} is defined by
	$$\fmaj(\pi):=r\cdot\maj(\pi)+\sum_{i=1}^n z_i.$$
\item $\rmaj$. The \emph{root major index}\cite{Haglund_Loehr_Remmel_05} is defined by
	$$\rmaj(\pi):=\maj(\pi)+\sum_{i=1}^n z_i\cdot \sigma_i.$$
\item $\fmaf$. The \emph{flag maf index}\cite{Faliharimalala_Zeng_11} is defined by
	$$\fmaf(\pi):=r\cdot \sum_{j=1}^k(i_j-j)+\fmaj(\tilde{\pi}),$$
	where $\{i_1,i_2,\ldots,i_k\}$ is the fix set $\Fix(\pi):=\{i:\pi_i=i\}$ and $\tilde{\pi}\in G(r,n-k)$ is obtained from $\pi$ by deleting fixed terms and then renumbering the remaining terms with order preserved.
\end{itemize}

For example, for $\pi=2^{[1]}\,1^{[3]}\,5\,4\,3^{[2]} \in G(4,5)$, we have $\fmaj(\pi)=38$, $\rmaj(\pi)=19$ and $\fmaf(\pi)=34$, with $\tilde{\pi} = 2^{[1]} 1^{[3]} 43^{[2]}$.

These three statistics have the same distribution over $G(r,n)$ and the following generating function (see \cite{Haglund_Loehr_Remmel_05, Faliharimalala_Zeng_11})
\begin{equation}\label{eq:Grn}
\sum_{\pi\in G(r,n)} q^{\fmaj(\pi)}=\sum_{\pi\in G(r,n)} q^{\fmaf(\pi)}
=\sum_{\pi\in G(r,n)} q^{\rmaj(\pi)}=\prod_{k=1}^n[rk]_q.
\end{equation}

Note that the statistic $\fmaj$ on $G(r,n)$ is a generalization of $\fmaj$ on $B_n$. 
However, it is no longer equidistributed with the length function for $r\geq 3$.

Now we define a new statistic \emph{root inversion index} $\rinv$ on $G(r,n)$.

\begin{defi} Given $\pi=(z,\sigma)\in G(r,n)$, the \emph{root inversion number} of $\pi$ is
$$\rinv(\pi):=\inv(\pi)+\sum_{i=1}^n z_i\cdot \sigma_i.$$
\end{defi}

For example, for $\pi=2^{[1]}\,1^{[3]}\,5\,4\,3^{[2]} \in G(4,5)$ we have $\rinv(\pi)=16$. 
Similar to the fact that $\inv$ and $\maj$ have a symmetric joint distribution
\begin{equation}\label{eq:invmaj_sym}
\sum_{\pi\in \mathfrak{S}_n}t^{\inv(\pi)}q^{\maj(\pi)}=\sum_{\pi\in \mathfrak{S}_n}q^{\inv(\pi)}t^{\maj(\pi)}
\end{equation}
over $\mathfrak{S}_n$~\cite{Foata_Sch_78}, we will prove that the statistics $\rinv$ and $\rmaj$ have a symmetric joint distribution over $G(r,n)$.

\begin{thm}\label{thm_invmaj}
The joint distribution of $\rmaj$ and $\rinv$ is symmetric over $G(r,n)$. 
That is,
$$\sum_{\pi\in G(r,n)}t^{\rinv(\pi)}q^{\rmaj(\pi)}= \sum_{\pi\in G(r,n)}q^{\rinv(\pi)}t^{\rmaj(\pi)}.$$
\end{thm}
\proof 
Recall that $Z(\pi)=\sum z_i$ for $\pi=(z,\sigma)\in G(r,n)$.
Let $\hat{Z}(\pi):=\sum z_i\cdot\sigma_i$ and therefore $\hat{Z}(\pi)=\hat{Z}(\tau)$ if $\pi=\tau\rho\in U_{r,n}\cdot \mathfrak{S}_n$.
\begin{eqnarray*}
&&\sum_{\pi\in G(r,n)}t^{\rinv(\pi)}q^{\rmaj(\pi)}
=\sum_{\pi\in G(r,n)}t^{\inv(\pi)+\hat{Z}(\pi)}q^{\maj(\pi)+\hat{Z}(\pi)}\\
&=&\sum_{\tau\in U_{r,n},\,\rho\in \mathfrak{S}_n}t^{\inv(\rho)+\hat{Z}(\tau)}q^{\maj(\rho)+\hat{Z}(\tau)}
=\sum_{\tau\in U_{r,n}}(tq)^{\hat{Z}(\tau)}\sum_{\rho\in \mathfrak{S}_n}t^{\inv(\rho)}q^{\maj(\rho)}\\
&\stackrel{(*)}{=}&\sum_{\tau\in U_{r,n}}(tq)^{\hat{Z}(\tau)}\sum_{\rho\in \mathfrak{S}_n}q^{\inv(\rho)}t^{\maj(\rho)}
=\sum_{\tau\in U_{r,n},\,\rho\in \mathfrak{S}_n}q^{\inv(\rho)+\hat{Z}(\tau)}t^{\maj(\rho)+\hat{Z}(\tau)}\\
&=&\sum_{\pi\in G(r,n)}q^{\inv(\pi)+\hat{Z}(\pi)}t^{\maj(\pi)+\hat{Z}(\pi)}
=\sum_{\pi\in G(r,n)}q^{\rinv(\pi)}t^{\rmaj(\pi)},
\end{eqnarray*}
where $(*)$ is because of the symmetric joint distribution of $\maj$ and $\inv$ over $\mathfrak{S}_n$. \qed

%%%%%%%%%%%%%%%%%%%%%%%%
\subsection{Signed counting polynomials}\label{Sec:FM_Results}
The following bivariate generating function, which is an extension of \eqref{eq:B(t,q)}, will play a key role in our discussion. 
The proof is again similar to that of Lemma~\ref{lem:r_bivariate} and is omitted.

\begin{lem}\label{lem:r&f} We have
$$R(t,q):=\sum_{\tau\in U_{r,n}}t^{Z(\tau)}q^{\hat{Z}(\tau)}=\prod_{k=1}^n[r]_{tq^k}.$$
\end{lem}

We are ready for the main results of this section. 
First we have $\rmaj$ v.s. $\fmaj$.

\begin{thm}\label{thm:fmajrmaj}
$$\sum_{\pi\in G(r,n)}t^{\rmaj(\pi)}q^{\fmaj(\pi)} =\prod_{k=1}^{n}[r]_{t^kq}[k]_{tq^r}.$$
\end{thm}
\proof 
The left-hand side is equal to
\begin{eqnarray*}
\sum_{\tau\in U_{r,n}}t^{\hat{Z}(\tau)}q^{Z(\tau)}\sum_{\rho\in \mathfrak{S}_n}t^{\maj(\rho)}q^{r\cdot\maj(\rho)},
\end{eqnarray*}
and the result follows by Lemma~\ref{lem:r&f} and \eqref{eq:A_Poincare}.
\qed

Then we have $\rinv$ v.s. $\rmaj$:
\begin{thm}\label{thm:rmajrinv}
$$\sum_{\pi\in G(r,n)}(-1)^{\rinv(\pi)}q^{\rmaj(\pi)} = \sum_{\pi\in G(r,n)}(-1)^{\rmaj(\pi)}q^{\rinv(\pi)} = \prod_{k=1}^n[r]_{(-q)^{k}}[k]_{(-1)^{k-1}q}.$$
\end{thm}
\proof
%Following the proof of Lemma \ref{thm_invmaj},
The left-hand side is equal to
\begin{eqnarray*}
\sum_{\tau\in U_{r,n}}(-q)^{\hat{Z}(\tau)}\sum_{\rho\in \mathfrak{S}_n}(-1)^{\inv(\rho)}q^{\maj(\rho)},
\end{eqnarray*}
and the result follows by Theorem~\ref{thm_invmaj}, Lemma~\ref{lem:r&f} and \eqref{eq:A_signMaho}.
\qed

And we have $\rinv$ v.s. $\fmaj$:

\begin{thm}\label{thm:fmajrinv}
$$\sum_{\pi\in G(r,n)}(-1)^{\rinv(\pi)}q^{\fmaj(\pi)} = \prod_{k=1}^n[r]_{(-1)^kq}[k]_{(-1)^{k-1}{q^r}}$$
and
$$\d\sum_{\pi\in G(r,n)}(-1)^{\fmaj(\pi)}q^{\rinv(\pi)} = \prod_{k=1}^n[r]_{-q^k}[k]_{(-1)^{(k-1)r}q}.$$
\end{thm}
\proof
By (\ref{eq:A_signMaho}) and Lemma \ref{lem:r&f},
\begin{eqnarray*}
&&\sum_{\pi\in G(r,n)}(-1)^{\rinv(\pi)}q^{\fmaj(\pi)}
=\sum_{\tau\in U_{r,n}}(-1)^{\hat{Z}(\tau)}q^{Z(\tau)}\sum_{\rho\in \mathfrak{S}_n}(-1)^{\inv(\rho)}q^{r\cdot\maj(\rho)}\\
&=& R(q, -1)\prod_{k=1}^{n}[k]_{(-1)^{k-1}{q^r}}
=\prod_{k=1}^n[r]_{(-1)^kq}[k]_{(-1)^{k-1}{q^r}}.
\end{eqnarray*}

Similarly,
\begin{eqnarray*}
&&\sum_{\pi\in G(r,n)}(-1)^{\fmaj(\pi)}q^{\rinv(\pi)}=
\sum_{\tau\in U_{r,n}}(-1)^{Z(\tau)}q^{\hat{Z}(\tau)}\sum_{\rho\in \mathfrak{S}_n}(-1)^{r\cdot\maj(\rho)}{q^{\inv(\rho)}}\\
&\stackrel{(*)}{=}& R(-1, q)\prod_{k=1}^{n}[k]_{(-1)^{(k-1)r}q}
=\prod_{k=1}^n[r]_{-q^k}[k]_{(-1)^{(k-1)r}q}.
\end{eqnarray*}
$(*)$: the second factor comes from \eqref{eq:A_Poincare} (or \eqref{eq:A_signMaho} and \eqref{eq:invmaj_sym}) if $r$ is even (or odd). \qed

The $\fmaf$ statistic is more elusive. In the following we will
prove signed counting results for even $r$ when
$(\stat_1,\stat_2)=(\fmaf,\rinv), (\fmaf, \rmaj),(\fmaf,\fmaj)$ or
$(\fmaj,\fmaf)$.

We need some preparations. 
For $\pi=(z,\sigma)\in G(r,n)$, denote the cardinality of $\Fix(\pi)$ by $\fix(\pi)$. 
Then $\fmaf(\pi)$ can be expressed as
\begin{align}\label{eq:maf_new}
\fmaf(\pi)&=r\cdot \left(\sum_{i\in \Fix(\pi)}i-{\fix(\pi)+1\choose 2}\right)+\fmaj(\tilde{\pi})\notag \\
&=r\cdot \left(\sum_{i\in\Fix(\pi)}i-{\fix(\pi)+1\choose 2}+\maj(\tilde{\pi})\right)+Z(\tilde{\pi}).
\end{align}

\begin{thm} \label{thm:fmafrmaj}
When $r$ is even, we have
$$\sum_{\pi\in G(r,n)}(-1)^{\fmaf(\pi)}q^{\rinv(\pi)}=
\sum_{\pi\in G(r,n)}(-1)^{\fmaf(\pi)}q^{\rmaj(\pi)}=\prod_{k=1}^{n}[r]_{-q^k}[k]_q.$$
\end{thm}
\proof
Let $(\Stat,\stat)$ denote $(\rinv,\inv)$ or $(\rmaj,\maj)$.
Since $r$ is even, by (\ref{eq:maf_new}), we have
\begin{align*}
\sum_{\pi\in G(r,n)}(-1)^{\fmaf(\pi)}q^{\Stat(\pi)}
=\sum_{\pi\in G(r,n)}(-1)^{Z(\tilde{\pi})}q^{\Stat(\pi)}
=\sum_{\tau\in U_{r,n}}(-1)^{Z(\tau)}q^{\hat{Z}(\tau)}\sum_{\rho\in \mathfrak{S}_n}q^{\stat(\rho)},
\end{align*}
%\begin{eqnarray*}
%&&\sum_{\pi\in G(r,n)}(-1)^{\fmaf(\pi)}q^{\Stat(\pi)}
%=\sum_{\pi\in G(r,n)}(-1)^{Z(\tilde{\pi})}q^{\Stat(\pi)}\\
%&=&\sum_{\tau\in U_{r,n}}(-1)^{Z(\tau)}q^{\hat{Z}(\tau)}\sum_{\rho\in \mathfrak{S}_n}q^{\stat(\rho)}
%=R(-1,q)\sum_{\rho\in \mathfrak{S}_n}q^{\stat(\rho)}=\prod_{k=1}^n[r]_{-q^k}[k]_q.
%\end{eqnarray*}
and the result follows by Lemma~\ref{lem:r&f} and \eqref{eq:A_Poincare}.
\qed

\begin{thm} \label{thm:fmaffmaj}
When $r$ is even, we have
$$\sum_{\pi\in G(r,n)}(-1)^{\fmaj(\pi)}q^{\fmaf(\pi)}
=\sum_{\pi\in G(r,n)}(-1)^{\fmaf(\pi)}q^{\fmaj(\pi)}=\prod_{k=1}^n[rk]_{-q}.$$
\end{thm}
\proof
Let $\stat_1$ be one of $\fmaj$ and $\fmaf$ and $\stat_2$ be the other.
Since $r$ is even, by (\ref{eq:maf_new}) and (\ref{eq:Grn}), we have
\begin{eqnarray*}
&&\sum_{\pi\in G(r,n)}(-1)^{\stat_1(\pi)}q^{\stat_2(\pi)} =\sum_{\pi\in G(r,n)}(-1)^{Z(\pi)}(-1)^{Z(\pi)}(-q)^{\stat_2(\pi)}\\
&=&\sum_{\pi\in G(r,n)}(-q)^{\stat_2(\pi)}=\prod_{k=1}^n[rk]_{-q}.
\end{eqnarray*}
Note that $Z(\tilde{\pi})=Z(\pi)$ by definition. \qed

Note that the results for $(\rmaj, \fmaf)$ and $(\rinv, \fmaf)$, as well as the odd $r$ cases in Theorem~\ref{thm:fmafrmaj} and Theorem~\ref{thm:fmaffmaj} are not known. 
The authors would like to know if these missing results can be completed.

%%%%%%%%%%%%%%%%%%%%%%%%%%%%%%%%%%%%%%%%%%%%%%%%%%%%%%%%%%%%%%%%%%%
%%%%%%%%%%%%%%%%%%%%%%%%%%%%%%%%%%%%%%%%%%%%%%%%%%%%%%%%%%%%%%%%%%%
\section{Concluding Notes}\label{sec:summary}
In this paper we proposed many new signed Mahonian polynomials for $B_n$, $D_n$ and $G(r,n)=G(r,1,n)$. 
In summary, for type $B_n$ we derived all signed Mahonian polynomials with respect to the statistics $\ell_B$, $\nmaj$, $\Fmaj$ and $\sor_B$, while for type $D_n$ we derived polynomials with respect to $\ell_D$, $\dmaj$ and $\sor_D$.
For $G(r,n)$ we first derived a closed formula for the 1-dim characters of $G(r,n)$ in terms of the length function, and as an application of it we obtained signed Mahonian polynomials with respect to $\ell$, $\lmaj$ and $\sor$.
This is novel in the sense that the better-known $\fmaj$ statistic is not equidistributed with the length function if $r\ge 3$ while $\lmaj$ is Mahonian. 

Caselli~\cite{Caselli_11} introduced the concept of projective reflection groups $G(r,p,s,n)$, which is a generalization of complex reflection groups $G(r,p,n)$, and Biagioli and Caselli~\cite{Biagioli_Caselli_12} derived signed polynomials with respect to the statistic $\fmaj$.
Hence a natural question is to see if one can have a Mahonian statistic $\lmaj$ or $\sor$ on $G(r,p,n)$ or $G(r,p,s,n)$ and have corresponding signed Mahonian polynomials.
%We leave these questions to the interested readers.

In the last section we considered the signed polynomials on $G(r,n)$ with respect to a pair of statistics with the same distribution $[r]_q[2r]_q\dots [nr]_q$, where the `sign' is decided by the parity of the first statistic. 
It would be interesting to complete the missing pieces in the last section, especially the $(\rmaj, \fmaf)$ and $(\rinv, \fmaf)$ cases. 
Also it is natural to investigate the bivariate generating functions, or even better, the multivariate generating functions with respect to these statistics. 
We leave these questions to the interested readers.

\bigskip

\noindent \textbf{Acknowledgements.}
The authors would like to express their gratitude to the referees for their valuable comments and suggestions on improving the presentation of this paper.

%%%%%%%%%%%%%%%%%%%%%%%%%%%%%%%%%%%%%%%%%%%%%%%%%%%%%%%%%%%%%%%%%%%%%%%%%%
%%%%%%%%%%%%%%%%%%%%%%%%%%%%%%%%%%%%%%%%%%%%%%%%%%%%%%%%%%%%%%%%%%%%%%%%%%
\rm

%%%%%%%%%%%%%%%%%%%%%%%%%%%%%%%%%%%%%%%%%%%%%%%%%%%%%%%%%%%%%%%%%%%%%%%%%%
%%%%%%%%%%%%%%%%%%%%%%%%%%%%%%%%%%%%%%%%%%%%%%%%%%%%%%%%%%%%%%%%%%%%%%%%%%

\bigskip
\begin{thebibliography}{24}

\bibitem{Adin_Brenti_Roichman_01}
R.M.~Adin, F.~Brenti, Y.~Roichman, \emph{Descent numbers and major indices for the hyperoctahedral group}, Special issue in honor of Dominique Foatas 65th birthday (Philadelphia, PA, 2000), Adv. Appl. Math. 27 (2001) 210--224.

\bibitem{Adin_Gessel_Roichman_05}
R.M.~Adin, I.~Gessel, Y.~Roichman, \emph{Signed Mahonians}, J. Combin. Theory Ser. A 109 (2005) 25--43.

\bibitem{Adin_Roichman_01}
R.M.~Adin, Y.~Roichman, \emph{The flag major index and group actions on polynomial rings}, European J. Combin. 22 (2001) 431--446.

\bibitem{Bagno_04}
E.~Bagno, \emph{Euler-Mahonian parameters on colored permutation groups}, S\'{e}minaire Lotharingien de Combinatoire 51 (2004), Article B51f.

\bibitem{Biagioli_03}
R.~Biagioli, \emph{Major and descent statistics for the even-signed permutation group}, Adv. in Appl. Math. 31 (2003) 163--179.

\bibitem{Biagioli_06}
R.~Biagioli, \emph{Signed Mahonian polynomials for classical Weyl groups}, European J. Combin. 27 (2006) 207--217.

\bibitem{Biagioli_Caselli_04}
R.~Biagioli and F.~Caselli, \emph{Invariant algebras and major indices for classical Weyl groups}, Proc. London Math. Soc. 88 (2004) 603--631.

\bibitem{Biagioli_Caselli_12}
R.~Biagioli and F.~Caselli, \emph{Weighted enumerations on projective reflection groups},
Adv. in Appl. Math. 48 (2012) 249--268.

\bibitem{Bjorner_Brenti_05}
A.~Bj\"{o}rner, F.~Brenti, \emph{Combinatorics Of Coxeter groups}, Springer (2005).

\bibitem{Caselli_11}
F.~Caselli, \emph{Projective reflection groups}, Israel J. Math. 185 (2011) 155--187.

\bibitem{ELW_15}
S.P.~Eu, Y.-H.~Lo and T.-L.~Wong, \emph{The sorting index on colored permutations and even-signed permutations}, Adv. in Appl. Math. 68 (2015) 18--50.

\bibitem{Faliharimalala_Zeng_11}
H.~Faliharimalala, J.~Zeng, \emph{Fix-Euler-Mahonian statistics on wreath products}, Adv. in Appl. Math. 46 (2011) 275--295.

\bibitem{Fire_08}
M.~Fire, \emph{Statistics on Wreath products}, November 2005, arXiv:math$\backslash$0409421v2.

\bibitem{Foata_Han_09}
D.~Foata, G.-H.~Han, \emph{New permutation coding and equidistribution of set-valued statistics}, Theoret. Comput. Sci., 410 (2009), 3743-3750.

\bibitem{Foata_Sch_78}
D.~Foata and M.P.~Sch\"{u}tzenberger, \emph{Major index and inversion number of permutations}, Math. Nachr. 83 (1978), 143--159.

\bibitem{Haglund_Loehr_Remmel_05}
J.~Haglund, N.~Loehr, J.B.~Remmel, \emph{Statistics on wreath products, perfect matchings and signed words}, European J. Combin. 26 (2005) 835--868.

\bibitem{Humphreys_90}
J.E.~Humphreys, \emph{Reflection Groups and Coxeter Groups}, Cambridge Studies in Advanced Mathematics, no. 29, Cambridge Univ. Press, Cambridge, 1990.

\bibitem{MacMahon_13}
P.A.~MacMahon, \emph{The indices of permutations and the derivation therefrom of functions of a single variable associated with the permutations of any assemblage of objects}, Amer. J. Math. 35 (1913) 281--322.

\bibitem{Petersen_11}
T.K.~Petersen, \emph{The sorting index}, Adv. in Appl. Math. 47 (2011) 615--630.

\bibitem{Reiner_95}
V.~Reiner, \emph{Descents and one-dimensional characters for classical Weyl groups}, Discrete Math. 140 (1995) 129--140.

\bibitem{Wachs_92}
M.~Wachs, \emph{An involution for signed Eulerian numbers}, Discrete Math. 99 (1992) 59-62.

\end{thebibliography}
\end{document}